\documentclass[reqno,11pt]{amsart}
\pdfoutput=1
\usepackage{amsmath}
\usepackage{amssymb}
\usepackage{graphicx}
\usepackage{amsfonts}
\usepackage[margin=0.93in, letterpaper ]{geometry}
\usepackage{lineno}
\usepackage[doublespacing]{setspace}
\usepackage{hyperref}

\hypersetup{
    colorlinks=true,
    citecolor=blue,
    linkcolor=blue,
    filecolor=blue,      
    urlcolor=blue,
}
\newtheorem{theorem}[equation]{Theorem}
\theoremstyle{plain}

\newtheorem{claim}[equation]{Claim}

\newtheorem{conjecture}[equation]{Conjecture}
\newtheorem{corollary}[equation]{Corollary}

\newtheorem{definition}[equation]{Definition}
\newtheorem{example}[equation]{Example}

\newtheorem{lemma}[equation]{Lemma}

\newtheorem{problem}[equation]{Problem}
\newtheorem{proposition}[equation]{Proposition}

\newtheorem{remark}[equation]{Remark}

\newtheorem{fact}[equation]{Fact}

\newtheorem{Observation}[equation]{Observation}
\numberwithin{equation}{subsection}
\input{tcilatex}

\begin{document}
\title[\textsf{Ranks for Representations\smallskip\ of }$GL_{n}$\textsf{\
Over Finite Fields, their Agreement, and Positivity}]{\textrm{Ranks for
Representations\smallskip\ of }$GL_{n}$\textrm{\ Over Finite Fields, their
Agreement, and Positivity of Fourier Transform}}
\author{\textsf{Shamgar Gurevich}}
\address{\textit{Department of Mathematics, University of Wisconsin,
Madison, WI 53706, USA.}}
\email{shamgar@math.wisc.edu}
\author{\textsf{Roger Howe}}
\address{\textit{Department of Mathematics, Yale University, New Haven, CT
06520, USA.}}
\email{roger.howe@yale.edu}

\begin{abstract}
In \cite{Frobenius1896}, it was shown that many important properties of a
finite group could be examined using formulas involving the \textit{%
character ratios }of group elements, i.e., the trace of the element acting
in a given irreducible representation, divided by the dimension of the
representation.

In \cite{Gurevich-Howe15} and \cite{Gurevich-Howe17}, the current authors
introduced the notion of \textit{rank} of an irreducible representation of a
finite classical group.

One of the motivations for studying rank was to clarify the nature of
character ratios for certain elements in these groups.

In fact in the above cited papers, two notions of rank were given. The first
is the Fourier theoretic based notion of $U$\textbf{-rank} of a
representation, which comes up when one looks at its restrictions to certain
abelian unipotent subgroups. The second is the more algebraic based notion
of \textbf{tensor rank} which comes up naturally when one attempts to equip
the representation ring of the group with a grading that reflects the
central role played by the few "smallest" possible representations of the
group.

In \cite{Gurevich-Howe17} we conjectured that the two notions of rank
mentioned just above agree on a suitable collection called "low rank"
representations.

In this note we review the development of the theory of rank for the case of
the general linear group $GL_{n}$ over a finite field $\mathbb{F}_{q}$, and
give a proof of the "agreement conjecture" that holds true for sufficiently
large $q$. Our proof is Fourier theoretic in nature, and uses a certain
curious positivity property of the Fourier transform of the set of matrices
of low enough fixed rank in the vector space of $m\times n$ matrices over $%
\mathbb{F}_{q}$.

In order to make the story we are trying to tell clear, we choose in this
note to follow a particular example that shows how one might apply the
theory of rank to certain counting problems.
\end{abstract}

\maketitle
\dedicatory{\smallskip\ \ \ \ \ \ \ \ \ \ \ \ \ \ \ \ \ \ \ \ \ \ \ \ \ \ \
\ \ \ \textrm{Dedicated to the memory of Tonny Springer}}

\section{\textbf{Introduction\label{S-In}}}

The Fourier theoretic study of a function on a finite abelian group via its
expansion as a linear combination of exponentials is by now a classical
example for the applications of harmonic analysis to pure and applied
mathematics \cite{Auslander-Tolimieri79}. This expansion has a well known
generalization to the study of class (i.e., invariant by conjugation)
functions on any finite group $G$. Indeed, let us denote by $\widehat{G}$
the set of (isomorphism classes of) complex finite dimensional irreducible
representations (\textit{irreps} for short) of $G$, and by%
\begin{equation}
\chi _{\pi },\text{ }\pi \in \widehat{G},  \label{IC}
\end{equation}%
the associated collection of irreducible characters, with $\chi _{\pi
}(g)=trace(\pi (g))$ for $g\in G$. In \ \cite{Schur1905} Schur formulated
his famous orthogonality relations for the collection (\ref{IC}), which
implies that it forms an orthogonal basis for the space of class functions
on $G$, equipped with the natural $G$-bi-invariant inner product on
functions on the group. This fact generalizes the abelian setting, and gives
birth to the theory of harmonic analysis on $G$, namely the investigation of
class functions on the group via their expansion as a linear combination of
irreducible characters.

As was already pointed out by Frobenius in \cite{Frobenius1896}, for many
interesting class functions on $G$ the expansion, as a linear combination of
irreducible characters, involves the normalized quantities%
\begin{equation}
\frac{\chi _{\pi }(g)}{\dim (\pi )},\text{ \ }\pi \in \widehat{G},\text{ }%
g\in G,  \label{CR}
\end{equation}%
called \textit{character ratios (CRs).}

So to make use of Frobenius's type formulas, it seems that one might benefit
from a solution to the following:\medskip

\textbf{Problem (Core problem of harmonic analysis on }$G$\textbf{). }%
Estimate the character ratios (\ref{CR}).\medskip

This note will focus on a particular example for the general and special
linear groups. We will show how to get precise information about the
character ratios for arbitrary representations of these groups for the
elements known as transvections. We introduce them now.

\subsection{\textbf{Example: Generation by Transvections\label{S-Example}}}

Consider the group $G=SL_{n}(\mathbb{F}_{q})$ of $n\times n$ matrices with
entries in a finite field $\mathbb{F}_{q}$ and determinant equal to one. For
this example let us assume that $n\geq 3.$ Inside $G$ we look at the
conjugacy class $C$ of the transvection\vspace{0.1cm}%
\begin{equation}
T\mathcal{=}%
\begin{pmatrix}
1 & 1 &  &  &  \\ 
& \ddots &  &  &  \\ 
&  & \ddots &  &  \\ 
&  &  & \ddots &  \\ 
&  &  &  & 1%
\end{pmatrix}%
,  \label{T}
\end{equation}%
\vspace{0.1cm}with $T_{ii}=1$ for $i=1,..,n;$ $T_{12}=1,$ and $T_{ij}=0$
elsewhere.\smallskip

It is not difficult to show (see \cite{Artin57}) that $C$ generates the
group $G$, and for a given element $g\in G$, we would like to understand in
how many ways it can be be obtained, i.e., for $\ell \geq 1$ what is the
cardinality of the set%
\begin{equation}
M_{\ell ,g}=\left\{ (c_{1},\ldots ,c_{\ell })\in C^{\ell };\text{ \ }%
c_{1}\cdot \ldots \cdot c_{\ell }=g\right\} \text{?}  \label{Mlg}
\end{equation}

Let us try to answer the above question for the "typical" elements of $G$.
Before we do so, let us recall some information regarding the conjugacy
class $C$.\medskip

\textbf{Facts. }The following hold\footnote{%
We write $\#(X)$ for the number of elements in a finite set $X$.}$^{\text{,}%
} $\footnote{%
The notation $a(q)=o(b(q))$ means that $a(q)/b(q)\rightarrow 0$ as $%
q\rightarrow \infty .$}$^{\text{,}}$\footnote{%
The notation $c(q)+o(...)$ stands for $c(q)+o(c(q)).$}\smallskip :

\begin{itemize}
\item \textit{Cardinality: }$\#(C)=q^{2n-2}+o(...)$ \cite{Artin57};\smallskip

\item \textit{Generation: \ }every element of $G$ can be written as a
product of no more than $n$ elements from $C$ \cite{Humphries80};

Moreover,

\item \textit{Most elements:} the "boundary" $\partial (G)$, of members of $%
G $ that one can't form by less than $n$ products from $C$, is 
\begin{equation*}
\partial (G)\text{ }=\text{ }\{g\in G;\text{ }\ker (g-I)\neq 0\}\text{,}
\end{equation*}%
in particular\footnote{%
We write $a(q)=O(b(q))$ if there is constant $A$ with $a(q)\leq A\cdot b(q)$
for all sufficiently large $q$.}, 
\begin{equation*}
\#(\partial (G))=\#(G)\cdot \left( \text{ }1-O(\frac{1}{q})\right) .
\end{equation*}
\end{itemize}

Of course, most of the elements of $\partial (G)$ are regular semi-simple,
i.e., are diagonalizable over some field extension of $\mathbb{F}_{q},$ and
have $n$ different eigenvalues there. These are our typical elements for the
example we are giving, and we would like to solve for them the following:

\begin{problem}[Generation]
\label{P-Gen}For a regular semi-simple element $g\in \partial (G)$, what is
the cardinality of $M_{\ell ,g}$ (\ref{Mlg})?
\end{problem}

\begin{figure}[h]\centering
\includegraphics
{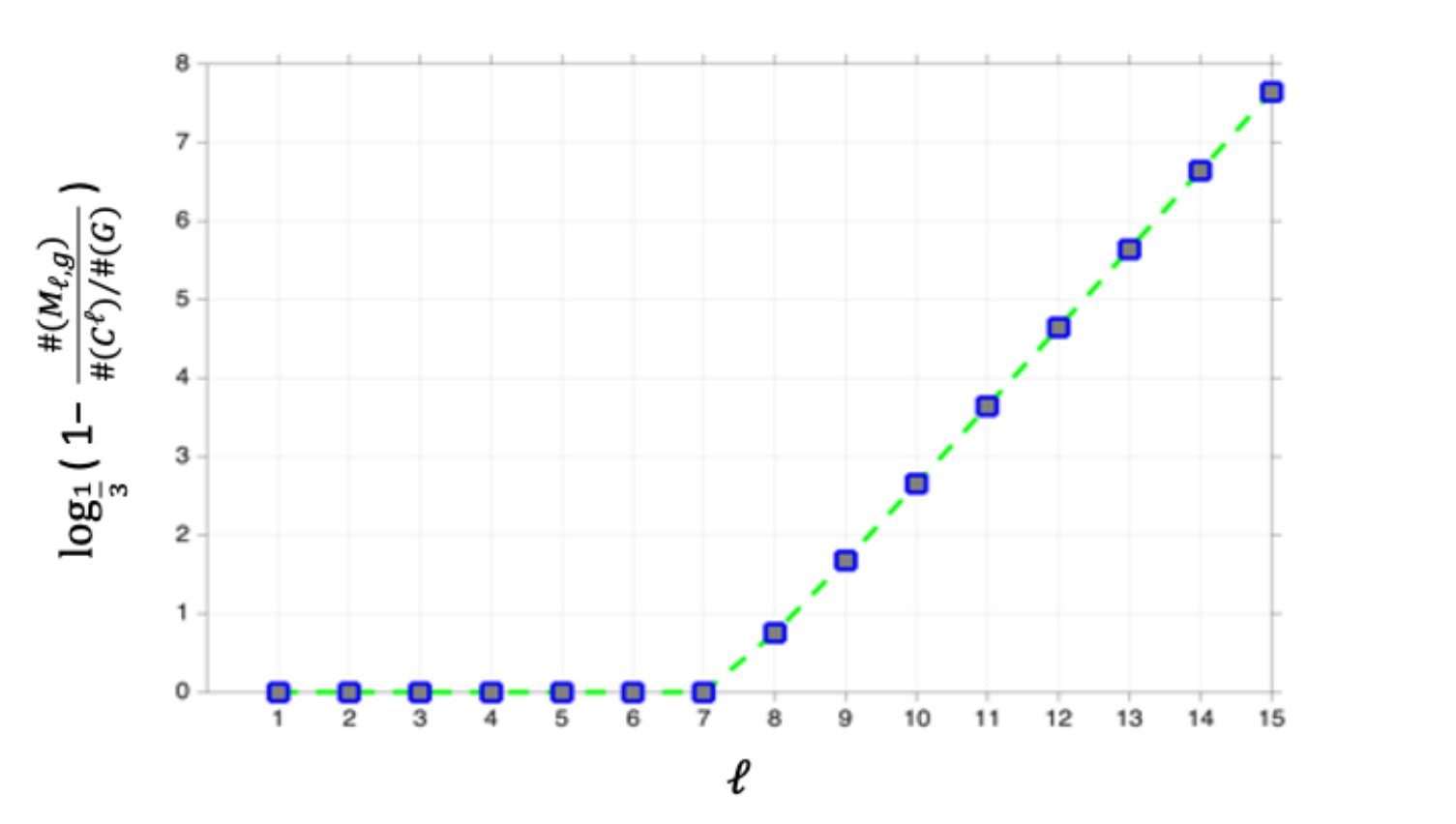}%
\caption{$\log _{\frac{1}{3}}$-scale of $(1-\frac{\#(M_{\ell ,g})}{%
\#(C^{\ell })/\#(G)})$ vs. $\ell $, for a typical element of $G=SL_{8}(%
\mathbb{F}_{3})$.}\label{gp-sl8-3}
\end{figure}

Note that, because we specialized to the case of a typical $g$, it makes
sense to expect (and probably not difficult to prove) that $\#(M_{\ell
,g})\rightarrow \#(C^{\ell })/\#(G)$, as $\ell \rightarrow \infty $. Before
we write down a precise statement, let us look at some numerics\footnote{%
The numerics in this note were done using the \href{http://magma.maths.usyd.edu.au/magma/}%
{Magma Computational Algebra System}.} for the ratio of $\#(M_{\ell ,g})$
and $\#(C^{\ell })/\#(G)$. Figure \ref{gp-sl8-3} illustrates, for the group $%
G=SL_{8}(\mathbb{F}_{3})$, how this quantity is close to being $1$ in $\log
_{\frac{1}{3}}$-scale. 
Let us elaborate a bit on what you see there. Of course, for $G=SL_{8}(%
\mathbb{F}_{3}),$ and our choice of $g,$ the set $M_{\ell ,g}$ is empty for $%
\ell <8$; but then the numerics shows that a "cutoff phenomenon" occurs,
namely, at the $\ell =8$ step the two quantities $\#(M_{\ell ,g})$ and $%
\#(C^{\ell })/\#(G)$ all of a sudden come close, and then at every
additional step they come closer by a multiple of $\frac{1}{3}$.

The numerical observations made just above, can be formulated and proved.
Indeed,

\begin{theorem}[Set-theoretic size]
\label{T-STS}\smallskip For a regular semi-simple $g\in \partial (G)$, we
have,\smallskip 
\begin{equation}
\#(M_{\ell ,g})=\frac{\#(C^{\ell })}{\#(G)}\cdot \left\{ 
\begin{array}{c}
1-O(\frac{1}{q}),\text{ \ \ \ \ \ \ \ \ \ \ \ \ \ \ \ \ if \ }\ell =n; \\ 
1-(\frac{2}{q}(\frac{1}{q})^{\ell -n}+o(...)),\text{ \ if \ }\ell >n.%
\end{array}%
\right.  \label{C-Mlg}
\end{equation}
\end{theorem}

\subsection{\textbf{A Geometric Analog of the Generation Problem}}

The set $M_{\ell ,g}$ (\ref{Mlg}) is in a natural way the set of $\mathbb{F}%
_{q}$-rational points of an algebraic variety $\mathbf{M}_{\ell ,g}$ defined
over $\mathbb{F}_{q}$, and both objects can teach us something about the
other (for basic terminology of algebraic geometry see \cite{Hartshorne77}).

A prototype example for the relation mentioned above---and relevant to our
story--- is given by the famous "Lang-Weil bound" \cite{Lang-Weil54}. It
relates the dimension of an (affine) algebraic variety $\mathbf{X}$ defined
over $\mathbb{F}_{q}$ (i.e., the set of solutions in an algebraic closure $%
\overline{\mathbb{F}}_{q}$ of $\mathbb{F}_{q}$ of a finite set of polynomial
equations over $\mathbb{F}_{q},$ and topology also given by polynomials over 
$\mathbb{F}_{q}$) and the cardinality of the set $X=\mathbf{X}(\mathbb{F}%
_{q})$ of its $\mathbb{F}_{q}$-rational points (i.e., the solutions in $%
\mathbb{F}_{q}$ of the polynomials defining $\mathbf{X}$). Here is a precise
formulation that will serve us well.

\begin{fact}[Lang-Weil bound]
\label{F-LW}The following are equivalent:\smallskip

\begin{enumerate}
\item \textit{Set-theoretic size: }$\#(X)=q^{d}+O(q^{d-\frac{1}{2}})$, for
some integer $d\geq 0$.\smallskip

\item \textit{Geometric size:} $\mathbf{X,}$ as a variety over $\overline{%
\mathbb{F}}_{q},$ has a unique irreducible component of maximal dimension $%
d=\dim (\mathbf{X)}$; all other components have smaller dimension.
\end{enumerate}
\end{fact}

In our case, we consider the algebraic group $\mathbf{G=SL_{n}}$ defined
over $\mathbb{F}_{q}$, and the conjugacy class $\mathbf{C\subset G}$ of the
transvection $T$ (\ref{T}). Then, for any $g\in \mathbf{G}$ we can form the
algebraic variety $\mathbf{M}_{\ell ,g}\subset \mathbf{C}^{\ell }$ in
exactly the same way as in (\ref{Mlg}). Moreover, $\mathbf{M}_{\ell ,g}$ is
defined over $\mathbb{F}_{q}$, and, indeed, $M_{\ell ,g}=\mathbf{M}_{\ell
,g}(\mathbb{F}_{q})$. So, in view of the Lang-Weil bound, a reasonable
geometric analog of Problem \ref{P-Gen}, might be the following:

\begin{problem}[Generation - geometric version]
\label{P-Gen-Geom}For regular semi-simple element $g\in \partial (\mathbf{G}%
) $ (i.e., no eigenvalue equal to $1$), compute the dimension of $\mathbf{M}%
_{\ell ,g}$, $\ell \geq n$, and the number of its irreducible components of
maximal dimension.
\end{problem}

To solve Problem \ref{P-Gen-Geom}, note that $\mathbf{M}_{\ell ,g}$ is the
fiber over $g$ of the multiplication morphism from $\mathbf{C}^{\ell }$ to $%
\mathbf{G}$. Hence, by the general "fiber dimension theorem" \cite%
{Hartshorne77}, for $\ell \geq n$, all components of $\mathbf{M}_{\ell ,g}$
have dimension $\geq \dim (\mathbf{C}^{\ell })-\dim (\mathbf{G})$. In
particular, looking on Fact \ref{F-LW}, we learn that Theorem \ref{T-STS}
implies the following:

\begin{corollary}[Geometric size]
\label{C-GS}Assume $g\in \partial (\mathbf{G})$ is regular semi-simple
element and $\ell \geq n$. Then $\mathbf{M}_{\ell ,g},$ as a variety over $%
\overline{\mathbb{F}}_{q},$ is irreducible of dimension 
\begin{eqnarray*}
\dim (\mathbf{M}_{\ell ,g}) &=&\dim (\mathbf{C}^{\ell })-\dim (\mathbf{G}) \\
&=&2\ell (n-1)-(n^{2}-1),
\end{eqnarray*}%
where in the last equality we used the fact (verified by a direct
computation) that $\dim (\mathbf{C})=2(n-1)$.
\end{corollary}

In fact, it will be interesting to find also a direct geometric proof of
Corollary \ref{C-GS}. However, it seems (compare the "error" term in Part
(1) of Fact \ref{F-LW} with the one appearing in Theorem \ref{T-STS}) that
the set-theoretic estimate we obtained is stronger than the geometric
information given in Corollary \ref{C-GS}.

It still left for us to explain why the counting statement appearing in
Theorem \ref{T-STS}, i.e., identity (\ref{C-Mlg}), is valid. For this, we
propose to use harmonic analysis.

\subsection{\textbf{Harmonic Analysis of the Generation Problem\label%
{S-HAofGP}}}

As a function of $g\in G=SL_{n}(\mathbb{F}_{q})$, the cardinality $%
\#(M_{\ell ,g})$, is a class function. The harmonic analytic expansion of
this function in irreducible characters can be computed explicitly. Indeed,

\begin{proposition}
\label{P-C-Mlg}We have,%
\begin{equation}
\#(M_{\ell ,g})=\frac{\#(C^{\ell })}{\#(G)}\cdot \left( 1+\underset{\mathbf{1%
}\neq \pi \in \widehat{G}}{\sum }\dim (\pi )\left( \frac{\chi _{\pi }(T)}{%
\dim (\pi )}\right) ^{\ell }\chi _{\pi }(g^{-1})\right) ,  \label{F-C-Mlg}
\end{equation}%
where $T$ is the transvection (\ref{T}).\smallskip
\end{proposition}

Formula (\ref{F-C-Mlg}) is well known \cite{Arad-Herzog-Stavi85,
Frobenius1896}; however, for the convenience of the reader we give another
verification in Appendix \ref{Pr-P-C-Mlg}.

The Formula (\ref{F-C-Mlg}), suggests proving Theorem \ref{T-STS} by
estimating the sum over the non-trivial representations,%
\begin{equation}
S_{\ell ,g}=\underset{\mathbf{1}\neq \pi \in \widehat{G}}{\sum }\dim (\pi
)\left( \frac{\chi _{\pi }(T)}{\dim (\pi )}\right) ^{\ell }\chi _{\pi
}(g^{-1}),  \label{N-T-Est}
\end{equation}%
and show that it is as small as the required "error" term in (\ref{C-Mlg}).

Recall that the element $g\in G$, appearing in the sum $S_{\ell ,g}$ (\ref%
{N-T-Est}), is regular semi-simple. For such generic elements the following
is known by \cite{Lusztig84} (and maybe can be deduced already from the work 
\cite{Green55}):

\begin{fact}
Suppose $g\in G$, is a regular semi-simple element. Then, there is a
constant $c,$ independent of $q$, such that for every irrep $\pi \in 
\widehat{G}$, we have,%
\begin{equation*}
\left\vert \chi _{\pi }(g)\right\vert \leq c\text{.}
\end{equation*}%
Moreover, one can take $c=n!$, the cardinality of the Weyl group $%
W_{G}=S_{n} $ of $G$.
\end{fact}

Looking back on (\ref{N-T-Est}) we see that, a possible approach for getting
the desired bound on $S_{\ell ,g}$ will be to have strong estimates on the
dimensions $\dim (\pi ),$ and, most importantly, on the character ratios $%
\frac{\chi _{\pi }(T)}{\dim (\pi )}$ of the irreps $\pi $ of $G=SL_{n}(%
\mathbb{F}_{q})$ at the transvection $T$ (\ref{T}).

In recent years d we have been developing a method that attempts to produce
this piece of information for the irreps of classical groups over finite
fields, and probably for character ratios of many other elements of interest.

\subsection{\textbf{Rank of a Representation}}

We want to estimate character ratio $\frac{\chi _{\pi }(T)}{\dim (\pi )},$
on the transvection $T$ (\ref{T}), for arbitrary irrep $\pi $ of $G=SL_{n}(%
\mathbb{F}_{q})$.

The group $G$ is a member of the family of reductive groups over finite and
local fields. The most popular method that people use to analyze
representations of such groups is the \textit{philosophy of cusp forms} \cite%
{Harish-Chandra70} put forward by Harish-Chandra in the 60s. In this
approach one studies the irreps of the group by means of certain basic
objects called \textit{cuspidal} representations\textit{. }It turns out that
cuspidality is a generic property, i.e., these irreps constitute a major
part of all irreps, and most of them are, in some sense, among the "largest".

The philosophy of cusp forms has had enormous success in establishing the
Plancherel formula for reductive groups over local fields \cite%
{Harish-Chandra84}, and leads to Lusztig's classification \cite{Lusztig84}
of the irreps of reductive groups over finite fields. However, analysis of
character ratios (CRs) seems to benefit from a different approach.

The attempts to estimate the CRs motivated us to introduce, in \cite%
{Gurevich-Howe15, Gurevich-Howe17}, a new way to think on the irreps of the
classical groups; a way in which the building blocks are the very few
"smallest" representations; in fact representations that may seem to be
anomalies in the cusp form approach.

In \cite{Gurevich-Howe15, Gurevich-Howe17} we explained that the choice of
looking on the irreps of a given classical group through the lens of its
smallest ones, reveals the existence of a pair of related invariants, which
we refer to by the label of "\textbf{rank}". Specifically, we have defined "$%
U$\textit{-rank" }and "\textit{tensor rank".} For example, for the case of
the group $GL_{n}(\mathbb{F}_{q})$, the $U$-rank of an irrep is an integer
between $0$ and $\frac{n}{2}$ or $\frac{n-1}{2}$, depending if $n$ is even
or odd, respectively, while tensor rank is an integer between $0$ and $n$.
In \cite{Gurevich-Howe17} we conjectured that in the case of $GL_{n}(\mathbb{%
F}_{q})$ for values in the range $0$ and $\frac{n}{2}-1$, for $n$ even, or $%
\frac{n-1}{2}-1$, for $n$ odd, these two invariants coincide.

\begin{remark}[General agreement conjecture and its meaning]
The conjectural agreement between $U$-rank and tensor rank for irreps of $%
GL_{n}$, is part of a general conjecture for all classical groups. Indeed,
in \cite{Gurevich-Howe17} we have defined these two invariants in the
mentioned generality, and conjectured that they agree on the collection of
"low" $U$--rank irreps. This would mean that tensor rank, which is defined
in a formal way using the representation (aka Grothendieck) ring, has a
concrete, down-to-earth meaning in terms of harmonic analysis on $G$ and its
subgroups. A future goal should be to extend this kind of interpretation to
representations of higher tensor rank.
\end{remark}

The first main goal of this note is to establish that agreement conjecture
for the group $GL_{n}(\mathbb{F}_{q}),$ for sufficiently large $q$. The
value of these two notions of rank is that they provide, in some sense, two
very different (to some extent complementary) reasons why certain analytic
properties of a representation, such as dimension an character ratio are, in
principle, what they are. So it is interesting and valuable to know that
these two notions in fact agree in the relevant range.\vspace{0.02in}

In particular, we show that,\vspace{0.02in}

\textbf{Theorem. }Fix $0\leq k\leq n$. Then for an irrep $\rho $ of $GL_{n}(%
\mathbb{F}_{q})$ of rank $k$, we have an estimate:\vspace{0.03in}%
\begin{equation}
\frac{\chi _{\rho }(T)}{\dim (\rho )}=\left\{ 
\begin{array}{c}
\text{ }\frac{1}{q^{k}}+o(...)\text{, \ \ \ \ if \ \ \ }k<\frac{n}{2};\text{
\ \ \ \ \ \ \ \ \ \ \ \ } \\ 
\text{\ \ \ \ \ \ } \\ 
\frac{c_{\rho }}{q^{k}}+o(...)\text{, \ \ \ \ if \ \ }\frac{n}{2}\leq k\leq
n-1;\text{\ \ } \\ 
\text{\ } \\ 
\text{ \ }\frac{-1}{q^{n-1}-1}\text{, \ \ \ \ \ \ \ \ if \ \ \ }k=n,\text{ \
\ \ \ \ \ \ \ \ \ \ \ }%
\end{array}%
\right. \vspace{0.03in}  \label{CRT}
\end{equation}
where $c_{\rho }$ is a certain integer (independent of $q$) combinatorially
associated with $\rho $.\vspace{0.05in}

\textbf{Remark. }For irreps $\rho $ of tensor rank $\frac{n}{2}\leq k\leq
n-1,$ the constant $c_{\rho }$ in (\ref{CRT}) might be equal to zero. In
this case, the estimate on $\frac{\chi _{\rho }(T)}{\dim (\rho )}$ is simply 
$o(\frac{1}{q^{k}}).$ However, it is typically non-zero, and in many cases
it is $1$.\smallskip

The estimates in (\ref{CRT}) induce similar results for the irreps of $%
G=SL_{n}(\mathbb{F}_{q})$. In particular, using some additional analytic
information, Theorem \ref{T-STS} on the cardinality of the set $M_{\ell ,g}$
(\ref{Mlg}) follows, and our introductory story is complete.\medskip

A first proof of the estimates (\ref{CRT}) appeared in \cite{Gurevich-Howe19}%
. However, the fact that irreps of the same rank have essentially the same
character ratio on the transvection (despite the fact that their dimensions
might differ by multiple of a large power of $q$) remains somewhat of a
surprise. In this note we clarify this phenomenon for low rank irreps (i.e.,
of rank $k<\frac{n}{2}$ or $\frac{n-1}{2},$ depending, respectively, if $n$
is odd or even). We give a clear picture why this is so, using the $U$-rank
realization of the notion of rank. This clarification is the second main
contribution of this note.\smallskip

\subsection{\textbf{Fourier Transform of Sets of Matrices of Fixed Rank}}

The third and final contribution of this note is an explicit formula for the
value of the Fourier transform of the set $(M_{n,n})_{k}$, of $n\times n$
matrices of rank $k\leq n,$ over a finite field $\mathbb{F}_{q}$, evaluated
at a rank one matrix $T$. The formula leads to an observation that for $k<n$%
, this value is \textbf{positive, }which turns out to be a significant
ingredient in our proof of the agreement conjecture mentioned
above.\smallskip

Let us write down the formula. Fixing an additive character $\psi \neq 1$ of 
$\mathbb{F}_{q}$, we have in the standard manner the associated Fourier
transform $f\mapsto \widehat{f}$ on the space of complex valued function on $%
M_{n,n}$, given by 
\begin{equation*}
\widehat{f}(B)=\sum_{A\in M_{n,n}}f(A)\psi (-trace(B^{t}\circ A)),\text{ }%
B\in M_{n,n}\text{.}
\end{equation*}%
Consider now the characteristic function $1_{(M_{n,n})_{k}}$ of the set $%
(M_{n,n})_{k}$. It is easy to see that its Fourier transform takes only real
values. Denote by $\Gamma _{n,k}$ the Grassmannian of all $k$-dimensional
subspaces of $\mathbb{F}_{q}^{n}$.\vspace{0.02in}

\textbf{Theorem. }The value of $\widehat{1}_{(M_{n,n})_{k}}$ on a rank one
matrix $T\in M_{n,n}$, is an integer that satisfies,%
\begin{eqnarray*}
\widehat{1}_{(M_{n,n})_{k}}(T) &=&\sum_{A\in (M_{n,n})_{k}}\psi
(trace(T^{t}\circ A)) \\
&=&(q^{2n-k}-2q^{n}+1)(\frac{\#(\Gamma _{n,k})^{2}\#(GL_{k})}{(q^{n}-1)^{2}}%
),
\end{eqnarray*}%
and in particular it is positive if $k<n$, and negative if $k=n$.\smallskip

After the acknowledgements and table of contents part, we proceed to the
body of the note, and start with a detailed discussion on character ratios
and the notion of $U$-rank for irreps of $GL_{n}(\mathbb{F}_{q})$.\medskip

\textbf{Acknowledgements. }The material presented in this note is based upon
work supported in part by the National Science Foundation under Grants No.
DMS-1804992 (S.G.) and DMS-1805004 (R.H.). \ 

We want to thank Steve Goldstein for the help with numerical aspects of the
project, part of which is reported in this note. Also, S.G. thanks Dima
Arinkin for several interesting discussions.

This note was written during 2020-21 partly while S.G. visited the school of
education at Texas A\&M University, and the math department at Yale
University, and he wants to thank these institutions, and especially Roger
Howe at TAMU and Yair Minsky at Yale.

\tableofcontents

\section{\textbf{Character Ratios and }$U$\textbf{-Rank\label{S-CRs-URank}}}

In Section \ref{S-Example} we described an example that motivated the need
to extract information on irreps of $SL_{n}(\mathbb{F}_{q})$. However, let
us start with a slightly better behaved group, namely the group $%
GL_{n}=GL_{n}(\mathbb{F}_{q})$ of $n\times n$ invertible matrices with
entries in $\mathbb{F}_{q}$. Moreover, for this group\footnote{%
In this note, for clarity, we denote irreps of $GL_{n}$ mostly by $\rho $
and of $SL_{n}$ mostly by $\pi .$} let us concentrate for a while only on
the problem of estimating the character ratios (CRs) on the transvection $T$
(\ref{T}),

\begin{equation}
\frac{\chi _{\rho }(T)}{\dim (\rho )},\ \ \rho \in \widehat{GL}_{n}.
\label{CR-T}
\end{equation}%
We want now to develop some intuition for how the quantity (\ref{CR-T})
behaves.

\subsection{\textbf{Numerics for Character Ratios vs. Dimension\label%
{S-N-CRvsD}}}

\begin{figure}[h]\centering
\includegraphics
{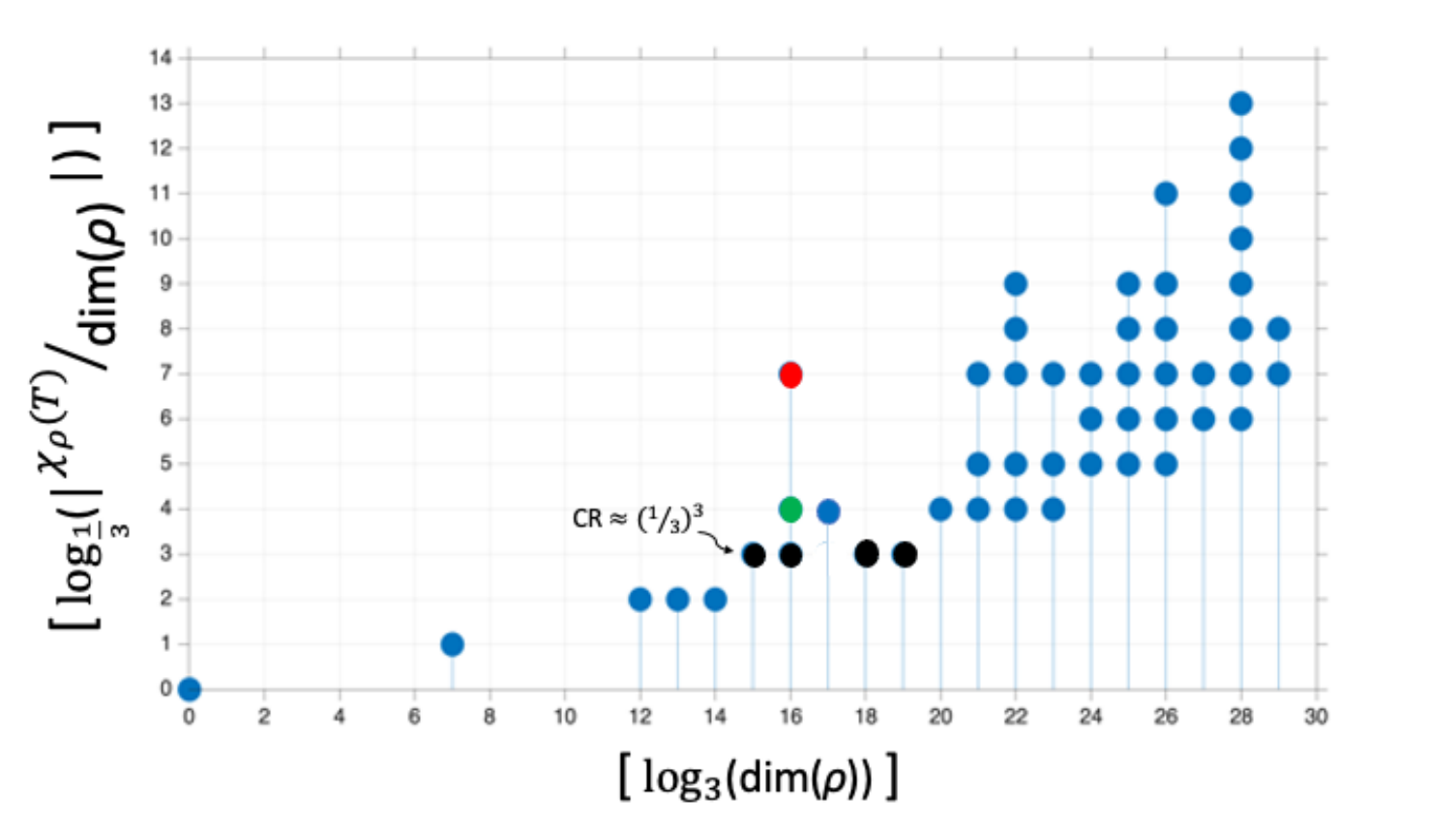}%
\caption{$\log _{\frac{1}{q}}\left\vert \frac{\protect\chi _{\protect\rho }(T%
\mathcal{)}}{\dim (\protect\rho )}\right\vert $ \ vs. $\ \log _{q}(\dim (%
\protect\rho ))$, $\ \protect\rho \in \protect\widehat{GL}_{8}(\mathbb{F}%
_{q}),$ $q=3$.}\label{cr-vs-dim-gl8}%
\end{figure}

Let us look at the numerics---appearing in Figure \ref{cr-vs-dim-gl8}---for
the group $GL_{8}(\mathbb{F}_{3})$. Let us explain a bit what 
appears there. For each irrep $\rho $ of this group we want to plot its CR
vs. its dimension. It can be deduced from \cite{Deligne-Lusztig76, Green55}
that, the dimensions of the irreps of groups like $GL_{n}(\mathbb{F}_{q})$
are certain polynomials in $q$, and as such have degrees. This integer can
be computed numerically, i.e., if $\rho \in \widehat{GL}_{8}(\mathbb{F}%
_{3}), $ then this degree is approximately (the nearest integer to) $\log
_{3}(\dim (\rho ))$, and this is what appears\footnote{%
We denote by $\left[ x\right] $ the nearest integer to the real number $x$.}
on the horizontal axis of Figure \ref{cr-vs-dim-gl8}. Next, looking closely
at the numerics for the CRs (\ref{CR-T}), one learns that they tend to come
in quantities which are powers of $\frac{1}{q}$, so it makes sense to plot
the (nearest integer of the) absolute value of them in $\log _{\frac{1}{3}}$%
-scale, and this is what appears on the vertical axis of Figure \ref%
{cr-vs-dim-gl8}.

What we can learn from Figure \ref{cr-vs-dim-gl8}? Let's read together part
of the data presented. The group $GL_{8}(\mathbb{F}_{q})$ has around $q$
irreps of dimension $1$, and of course their character ratios are of size
around $1$, this is the blue dot at $(0,0)$ in the figure. After that, we
have the around $q^{2}$ (more details later on why this is the cardinality)
irreps of dimension around $q^{7}$, all of them seem to have CR $\approx 
\frac{1}{q}$. Next, we have around $q^{3}$ irreps of dimensions (already
some variation) from around $q^{12}$ to $q^{14}$, but CR "exactly" $\frac{1}{%
q^{2}}$. Let us read one more layer, we have the black circles in Figure \ref%
{cr-vs-dim-gl8}, of around $q^{4}$ irreps, and dimensions vary (by multiple
of a quite large power of $q$) from around $q^{15}$ to $q^{19}$, but the CRs
of these irreps are nearly the same, of order of magnitude $\frac{1}{q^{3}}$%
. Finally, another look at Figure \ref{cr-vs-dim-gl8} reveals a collection
of irreps (see the black, green, and red circles above $16$ there) all of
them have the same dimension, around $q^{16},$ but for some reason they have
very different CRs (respectively, $\frac{1}{q^{3}},$ $\frac{1}{q^{4}},$ and $%
\frac{1}{q^{7}}$).

In summary, based on the numerics appearing in Figure \ref{cr-vs-dim-gl8},
we can make the following reasonable:

\begin{Observation}
In general, the character ratios (\ref{CR-T}) are not strictly controlled by
the dimensions of the irreps.
\end{Observation}

In the literature we are aware of, the CRs are estimated using information
on the dimensions of the irreps (for example see \cite%
{Bezrukavnikov-Liebec-Shalev-Tiep18}); so in general (e.g., on elements like
the transvection $T$ (\ref{T})) these estimates cannot be optimal.

In \cite{Gurevich-Howe15, Gurevich-Howe17} we initiated the study of a pair
of related invariants that seem to do better job than dimension in
controlling the CRs (\ref{CR-T}). We proceed to discuss the first of these
invariants.

\subsection{$U$\textbf{-Rank: Motivation, Intuition, and Formal Definition 
\label{S-UR}}}

Look again on Figure \ref{cr-vs-dim-gl8}. What makes the irreps of $GL_{8}(%
\mathbb{F}_{3})$ with character ratio of order of magnitude $\frac{1}{3^{k}}$%
, for $k<4$, a family?, i.e., what puts them together?

More generally, we want an invariant that in some sense "knows" which irrep $%
\rho $ of $GL_{n}=GL_{n}(\mathbb{F}_{q})$ has%
\begin{equation*}
\left\vert \frac{\chi _{\rho }(T)}{\dim (\rho )}\right\vert =\frac{1}{q^{k}}%
+o(...),\ \ \text{for}\ \ k<\frac{n-1}{2},
\end{equation*}

We want to give some intuition for the invariant we propose for such a job.
Let us start with some data (for basic notions from the theory of algebraic
groups see \cite{Borel69}). Consider the vector space $\mathbb{F}_{q}^{n},$
and the "first $m$-coordinates" subspace 
\begin{equation}
\mathbb{F}_{q}^{n}\supset X_{m}=\left\{ 
\begin{pmatrix}
x_{1} \\ 
\vdots \\ 
x_{m} \\ 
0 \\ 
\vdots \\ 
0%
\end{pmatrix}%
;\text{ }x_{j}\in \mathbb{F}_{q}\right\} .  \label{Xm}
\end{equation}%
This allows to define three subgroups of $GL_{n}$, that play an important
role in our story. The first is the stabilizer of $X_{m}$, i.e., the
parabolic subgroup%
\begin{equation}
P_{m}=Stab_{GL_{n}}(X_{m})=\{g\in GL_{n};\text{ }g(X_{m})=X_{m}\},
\label{Pm}
\end{equation}%
of elements of $GL_{n}$, that take the subspace $X_{m}$ to itself. Note that,

\begin{equation*}
P_{m}=\left\{ 
\begin{pmatrix}
C & A \\ 
0 & D%
\end{pmatrix}%
\right\} ,
\end{equation*}%
where $A\in M_{m,(n-m)}$---the space of $m\times (n-m)$ matrices---and $C\in
GL_{m},$ $D\in GL_{n-m}$.

In particular, we have a (split) short exact sequence of groups (with
obvious maps): 
\begin{equation}
1\rightarrow \overset{U_{m}}{\overbrace{\left\{ 
\begin{pmatrix}
I_{m} & A \\ 
0 & I_{n-m}%
\end{pmatrix}%
\right\} }}\hookrightarrow \overset{P_{m}}{\overbrace{\left\{ 
\begin{pmatrix}
C & A \\ 
0 & D%
\end{pmatrix}%
\right\} }}\twoheadrightarrow \overset{L_{m}}{\overbrace{\left\{ 
\begin{pmatrix}
C & 0 \\ 
0 & D%
\end{pmatrix}%
\right\} }}\rightarrow 1,  \label{Levi-decomp}
\end{equation}%
where $U_{m}$ and $L_{m}$, are called, respectively, the \textit{unipotent
radical} and \textit{Levi component} of $P_{m}=U_{m}\cdot L_{m}$. Here $%
I_{m} $ and $I_{n-m},$ are the identity matrices of order $m$ and $n-m$,
respectively.

Note that the group $U_{m}$ is commutative and naturally isomorphic to the
vector space $M_{m,(n-m)}$, using the map 
\begin{equation*}
\begin{pmatrix}
I_{m} & A \\ 
0 & I_{n-m}%
\end{pmatrix}%
\longmapsto \text{ }A.
\end{equation*}%
In particular, we might, and in many cases will, think of elements of $U_{m}$
as $m\times (n-m)$ matrices, and write $A\in U_{m}$, for a matrix $A\in
M_{m,(n-m)}$.

Next, let us denote by $\widehat{U}_{m}$ the \textit{Pontryagin dual} of the
commutative group $U_{m}$, consisting of all of its characters (one-dim
reps). Then, in the standard manner, fixing a non-trivial additive character 
$1\neq \psi $ of the field $\mathbb{F}_{q}$, one gets an isomorphism%
\begin{equation*}
\left\{ 
\begin{array}{c}
U_{m}\widetilde{\longrightarrow }\widehat{U}_{m}; \\ 
B\mapsto \psi _{B},%
\end{array}%
\right.
\end{equation*}%
where for $B\in U_{m}$, with transpose $B^{t}$, we define $\psi _{B}(A)=\psi
(trace(B^{t}A)),$ $A\in U_{m}$.

In particular, we have a notion of rank for every representation of $U_{m}$
as follows:

\begin{definition}
We define,

\begin{enumerate}
\item the \textit{rank of a character }$\psi _{B}\in \widehat{U}_{m},$ to be 
$rank(B)$;\smallskip

and,

\item the \textit{rank of a representation} $\varrho $ of $U_{m}$, to be the
maximum over the ranks of characters that appear in $\varrho $.
\end{enumerate}
\end{definition}

For simplicity of exposition, let us now restrict (however, see Remark \ref%
{R-U-rank} below) our attention to the case, 
\begin{equation}
U=U_{\left\lfloor \frac{n}{2}\right\rfloor }\text{, \ with }\left\lfloor
x\right\rfloor =\text{the largest integer below }x\text{.}  \label{U}
\end{equation}

\begin{example}[Fourier transform of rank $k$ matrices]
\label{E-FT} Fix integer $0\leq k\leq \left\lfloor \frac{n}{2}\right\rfloor
, $ and consider the collection 
\begin{equation*}
\mathcal{O}_{k}\subset U,
\end{equation*}%
of rank $k$ matrices in $U$, and by 
\begin{equation}
\varrho _{\mathcal{\mathcal{O}}_{k}}=\dsum\limits_{B\in \mathcal{O}_{k}}\psi
_{B},  \label{rho-Or}
\end{equation}%
the representation of $U$, which is the direct sum of all characters
corresponding to the members of $\mathcal{O}_{k}$. Then,

\begin{itemize}
\item $rank(\varrho _{\mathcal{\mathcal{O}}_{k}})=k.$

\item $\dim (\varrho _{\mathcal{\mathcal{O}}_{k}})=\#(\mathcal{O}_{k}).$

\item The character $\chi _{\varrho _{\mathcal{\mathcal{O}}_{k}}}$ of $%
\varrho _{\mathcal{\mathcal{O}}_{k}}$, satisfies, as a function on $U,$%
\begin{equation}
\chi _{\varrho _{\mathcal{\mathcal{O}}_{k}}}=\dsum\limits_{B\in \mathcal{O}%
_{k}}\psi _{B}=\widehat{1}_{\mathcal{O}_{k}},  \label{Chi_rho_Or}
\end{equation}%
where $\widehat{1}_{\mathcal{O}_{k}}$ is the Fourier transform (with respect
to the, previously fixed, additive character $\psi $ of $\mathbb{F}_{q}$) of
the characteristic function of $\mathcal{O}_{k}$.

\item The transvection $T$ (\ref{T}) has a $GL_{n}$-conjugate in $U,$ that
for simplicity we will also denote by $T$. Later, in Appendix \ref{A-FT}, we
will show that, the value of the character sum (\ref{Chi_rho_Or}) at $T$ is:

\begin{enumerate}
\item positive, if $k<\left\lfloor \frac{n}{2}\right\rfloor ;$ in fact, we
obtain an explicit formula for $\widehat{1}_{\mathcal{O}_{k}}(T)$ that
implies, 
\begin{equation}
\frac{\widehat{1}_{\mathcal{O}_{k}}(T)}{\#(\mathcal{O}_{k})}=\frac{1}{q^{k}}%
+o(...).  \label{CR-Or}
\end{equation}%
\smallskip and,

\item negative, if $k=\left\lfloor \frac{n}{2}\right\rfloor ;$ in fact we
get, 
\begin{equation*}
\frac{\widehat{1}_{\mathcal{O}_{\left\lfloor \frac{n}{2}\right\rfloor }}(T)}{%
\#(\mathcal{O}_{\left\lfloor \frac{n}{2}\right\rfloor })}=\frac{-1}{%
q^{\left\lfloor \frac{n}{2}\right\rfloor }}+o(...).
\end{equation*}
\end{enumerate}
\end{itemize}
\end{example}

In particular, we see that if $\rho $ is a rep of $GL_{n}$ which is
"supported" just on one orbit $\mathcal{O}_{k}\subset U$, i.e., $\rho $
restricted to $U$ satisfies $\rho _{|U}=m_{k}\cdot \varrho _{\mathcal{%
\mathcal{O}}_{k}}$ for some integer $m_{k}>0,$ then,\smallskip

(a) $\dim (\rho )=m_{k}\cdot \#(\mathcal{O}_{k});\smallskip \smallskip $

and, using (\ref{CR-Or}),\smallskip\ for $k<\left\lfloor \frac{n}{2}%
\right\rfloor ,$

(b) $\frac{\chi _{\rho }(T\mathcal{)}}{\dim (\rho )}=\frac{\widehat{1}_{%
\mathcal{O}_{k}}(T)}{\#(\mathcal{O}_{k})}=\frac{1}{q^{k}}+o(...).\smallskip $

The discussion above suggests an invariant that might explain the behavior
of the character ratio at the transvection. Indeed, take a rep $\rho $ of $%
GL_{n}$, and look at its restriction $\rho _{|U}$ to $U$. This has the
following description:

\begin{proposition}
\label{P-mult-rho_U}Characters of $U$of the same rank, appear in $\rho _{|U}$
with the same multiplicity, i.e., 
\begin{equation}
\rho _{|U}=\dsum\limits_{r=0}^{\left\lfloor \frac{n}{2}\right\rfloor
}m_{r}\cdot \varrho _{\mathcal{\mathcal{O}}_{r}},  \label{rho_U}
\end{equation}%
for some non-negative integers $m_{r}$, where $\varrho _{\mathcal{\mathcal{O}%
}_{r}}$ is given by Formula (\ref{rho-Or}).
\end{proposition}

For a proof of Proposition \ref{P-mult-rho_U} see Appendix \ref%
{P-P-mult-rho_U}.\smallskip

So, motivated by (b) in Example \ref{E-FT} above, we introduce the key
notion:

\begin{definition}
\label{D-URank}The \underline{$U$\textbf{-rank}} of a representation $\rho $
of $GL_{n}$, is the maximal $k$, $0\leq k\leq \left\lfloor \frac{n}{2}%
\right\rfloor $, such that $m_{k}\neq 0$, in (\ref{rho_U}).
\end{definition}

We will write $U$-$rank(\rho )=k$, or $rank_{U}(\rho )=k,$ to denote that a
rep $\rho $ of $GL_{n}$, has $U$-rank $k,$ and will use the notation $(%
\widehat{GL}_{n})_{U,k}$, for the set of all irreps with this
property.\smallskip

Sometimes we will call representations of $U$-rank less than $\left\lfloor 
\frac{n}{2}\right\rfloor ,$ \textit{low }$U$-\textit{rank representations.}

\begin{remark}
\label{R-U-rank}It was shown in \cite{Gurevich-Howe17} that for a low $U$%
-rank representation $\rho $ of $GL_{n}$, the value $rank_{U}(\rho )$ is
"independent" of $U,$ in the following sense. Consider in $GL_{n}$ a general
parabolic subgroup of block upper triangular matrices:%
\begin{equation*}
\left\{ 
\begin{pmatrix}
C_{1} & A_{12} & A_{13} & \cdots & \cdots & A_{1l} \\ 
& C_{2} & A_{23} & \cdots & \cdots & A_{2l} \\ 
&  & C_{3} & \cdots & \cdots & A_{3l} \\ 
&  &  & \cdots & \cdots & \cdots \\ 
&  &  &  & \cdots & \cdots \\ 
&  &  &  &  & C_{l}%
\end{pmatrix}%
\right\} .
\end{equation*}%
If the matrices $C_{j}$ are of size $m_{j}\times m_{j}$, then for a fixed $%
1\leq i<j\leq l$, the collection of matrices $\{A_{ij}\}$ forms a
(unipotent) subgroup of $GL_{n}$ isomorphic to $m_{i}\times m_{j}$ matrices.
Let us call such subgroups of $GL_{n}$\textit{\ standard matrix subgroups}.
The unipotent radicals $U_{m},$ introduced above, are examples of \ such
subgroups, and in particular the group $U$ (\ref{U}). The point is that you
can develop the theory discussed above using the restrictions to each of the
standard matrix subgroups. In particular, we have an induced notion of rank
for representations of $GL_{n}$ which is associated with each of these
subgroups, and a definition of what does it means for a representation to be
of low rank in each case. In \cite{Gurevich-Howe17} we showed that if we
have a representation $\rho $ of $GL_{n}$, and two standard matrix subgroups 
$U_{l,m}\simeq M_{l,m},$ $U_{l^{\prime },m^{\prime }}\simeq M_{l^{\prime
\prime },m^{\prime }}\subset GL_{n}$, such that $\rho $ is of low $U_{l,m}$%
-rank $k,$ i.e., $rank_{U_{l,m}}(\rho )=k<l,m$, and $U_{l^{\prime
},m^{\prime }}$ is big enough, i.e., $l\leq l^{\prime }$ and $m\leq
m^{\prime }$, then 
\begin{equation*}
rank_{U_{l,m}}(\rho )=rank_{U_{l^{\prime },m^{\prime }}}(\rho )\text{.}
\end{equation*}%
In particular, for a $U$-rank $k<\left\lfloor \frac{n}{2}\right\rfloor $
irrep, any large enough standard matrix subgroup can be used to detect this
invariant. The meaning of this might be that, from the point of view of
standard matrix subgroups the notion of rank developed above seems somehow
canonical.
\end{remark}

Let us look at some numerics that provide further evidence that we are on
the right track.

\subsection{\textbf{Numerics for Character Ratios vs. }$U$-Rank}

\begin{figure}[h]\centering
\includegraphics
{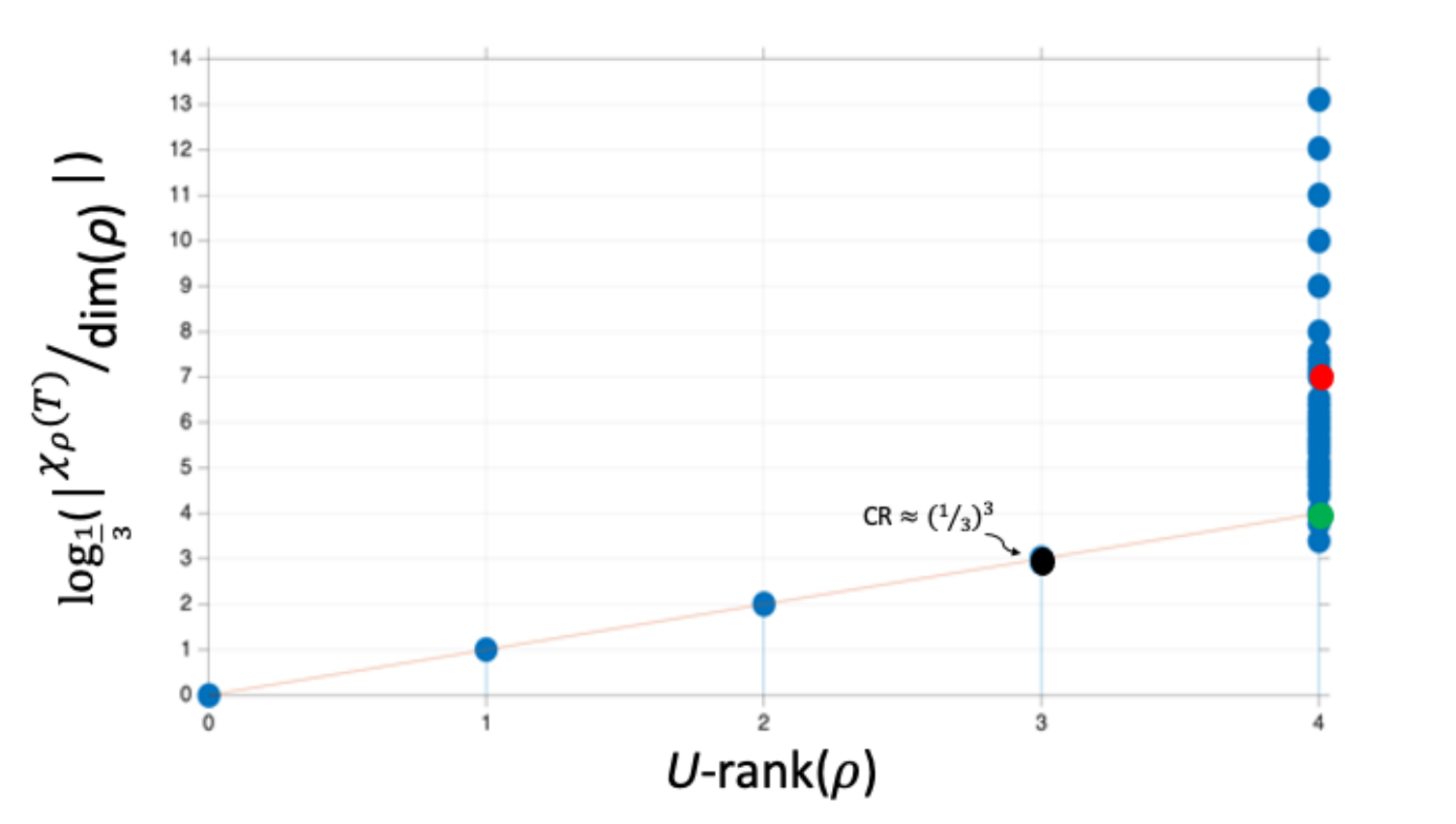}%
\caption{$\log _{\frac{1}{q}}$-scale of CRs\ vs. $U$-rank for $\protect\rho %
\in \protect\widehat{GL}_{8}(\mathbb{F}_{q}),$ $q=3.$}\label{cr-vs-ur-gl8}%
\end{figure}

As we saw in Part (b) of Example \ref{E-FT} above, if $\rho $ is an irrep of 
$GL_{n}$, for which its restriction to $U$, is supported solely on the rank $%
k$ matrices in $U$, for $k<$ $\left\lfloor \frac{n}{2}\right\rfloor $, then
its character ratio at the transvection (\ref{T}) is around $\frac{1}{q^{k}}$%
. This leads us to define the notion of $U$-rank $k$, with the hope that if $%
\rho $ is a representation of that rank, then the contribution of the lower
orbits in $U,$ on which it might also be supported, will not contribute much
to its CR on the transvection. The numerics done for $GL_{8}(\mathbb{F}_{3})$
and appear in Figure \ref{cr-vs-ur-gl8}, illustrate the fact that this is
probably true for all (and probably only for) the irreps of rank $k$, as
long as $k<$ $\left\lfloor \frac{n}{2}\right\rfloor $. In particular, Figure %
\ref{cr-vs-ur-gl8} suggests that 
for the irreps of $GL_{n}$ of rank less than $\left\lfloor \frac{n}{2}%
\right\rfloor $, the $U$-rank invariant does a better job than (compare with
Figure \ref{cr-vs-dim-gl8}) dimension in controlling the CRs at the
transvection. In that range it puts (compare the black circles in Figures %
\ref{cr-vs-ur-gl8} and \ref{cr-vs-dim-gl8}) the irreps of the "same" CRs
together.

The above numerical observations can be quantified precisely and proved.
This is part of what we do next.

\section{\textbf{Analytic Information on }$U$-\textbf{Rank }$k$\textbf{\
Irreps of }$GL_{n}$\textbf{\label{S-AI-UR}}}

To bound the sum (\ref{N-T-Est})---discussed in Section \ref{S-HAofGP}---we
would like to formulate statements on certain analytic properties of irreps
of $GL_{n}$. In particular, on their character ratios at the transvection $T$
(\ref{T}), and on their dimensions. We give now precise information on these
quantities---in fact sharp estimates in term of the rank---for the irreps of 
$U$-rank $k,$ with $k<$ $\left\lfloor \frac{n}{2}\right\rfloor $. In
addition, we calculate the number of such irreps.

The results follow from the analog results for tensor rank $k$ irreps
formulated in Theorem \ref{T-AI-GLn}, and the fact, given in Theorem \ref%
{T-Agreement}, that if $k<\left\lfloor \frac{n}{2}\right\rfloor ,$ then for
an irrep of $GL_{n}$ being of $U$-rank $k$ is the same thing as being of
tensor rank $k$.

\subsection{\textbf{Character Ratios on the Transvection}}

For this quantity we have,

\begin{theorem}
\label{T-CRs-URank}Suppose $k<\left\lfloor \frac{n}{2}\right\rfloor $, and $%
\rho $ is an irrep of $GL_{n}$ of $U$-rank $k$. Then,%
\begin{equation}
\frac{\chi _{\rho }(T)}{\dim (\rho )}=\ \frac{1}{q^{k}}+o(...).
\label{CRs-T-U-rank}
\end{equation}
\end{theorem}

Note that (\ref{CRs-T-U-rank}) is a formal validation to some of the
phenomena that Figure \ref{cr-vs-ur-gl8} illustrates.

\subsection{\textbf{Dimensions of Irreps\label{S-Dim-Urank-k}}}

As we already remarked earlier, although the CRs (\ref{CRs-T-U-rank}) of the
irreps of $GL_{n}$ of $U$-rank $k,$ $k<\left\lfloor \frac{n}{2}\right\rfloor 
$, are approximately the same, their dimensions might vary by multiple of a
large power of $q$. 
\begin{figure}[h]\centering
\includegraphics
{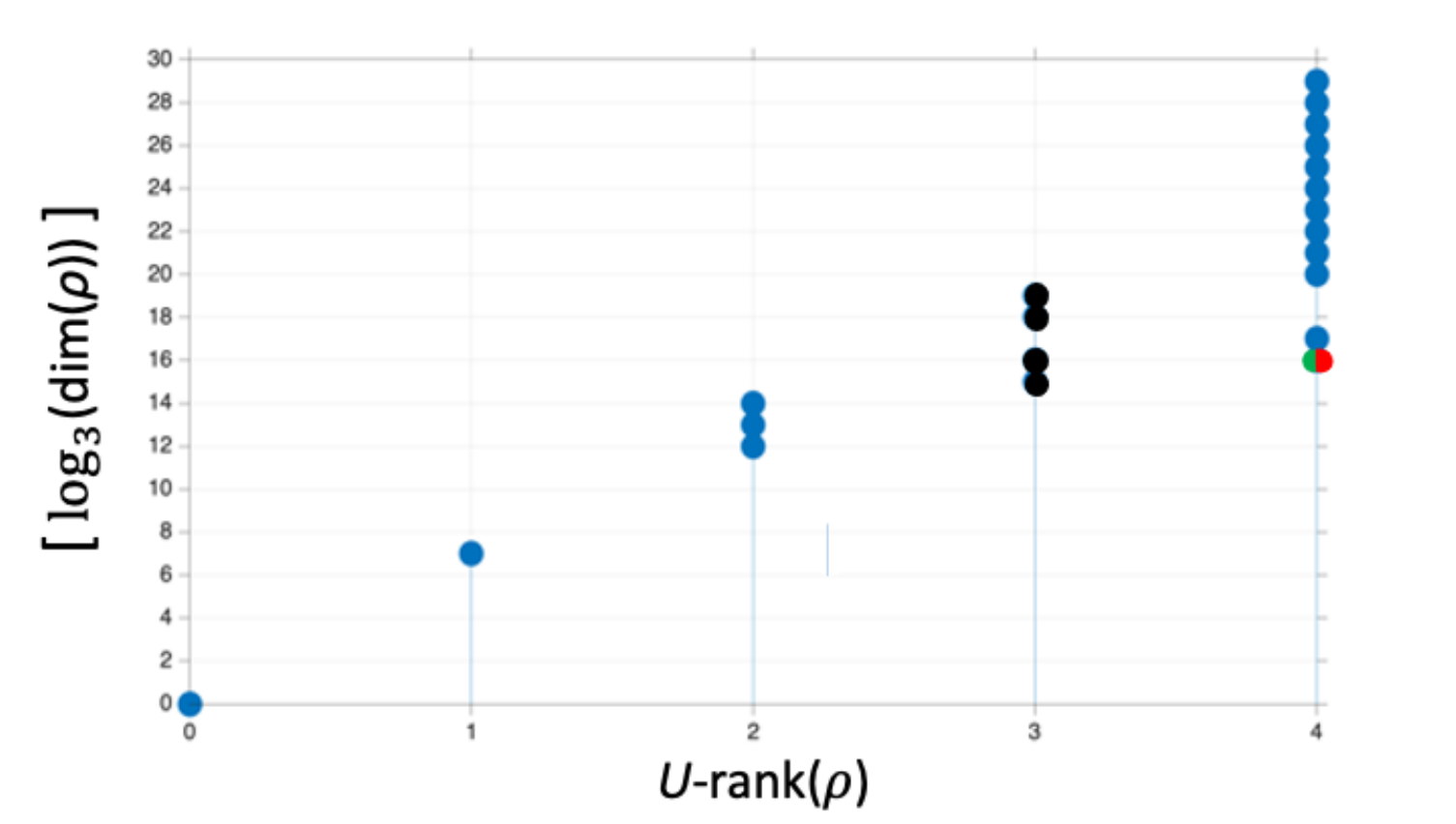}%
\caption{$\log _{q}$-scale of dimension vs. $U$-rank for irreps $\protect%
\rho $ of $GL_{8}(\mathbb{F}_{q}),$ $q=3.$}\label{dim-vs-ur-gl8}%
\end{figure}%

Figure \ref{dim-vs-ur-gl8} illustrates (compare with Figures \ref%
{cr-vs-ur-gl8} and \ref{cr-vs-dim-gl8}) the distribution of the dimensions
of the irreps of $GL_{8}(\mathbb{F}_{3})$ within each given $U$-rank $k$.
What you see there, can be formulated and proved in general for low rank
irreps of $GL_{n}$. Indeed, we have the following sharp lower and upper
bounds in term of the $U$-rank.

\begin{theorem}
\label{T-Dim-U-rank} Suppose $k<\left\lfloor \frac{n}{2}\right\rfloor $, and 
$\rho $ is an irrep of $GL_{n}$ of $U$-rank $k$. Then,%
\begin{equation}
q^{k(n-k)}+o(...)\leq \dim (\rho )\leq q^{k(n-k)+\frac{k(k-1)}{2}%
}+o(...).\medskip \medskip  \label{Dim-U-rank}
\end{equation}%
Moreover, the upper and lower bounds in (\ref{Dim-U-rank}) are attained.
\end{theorem}

Looking on the upper and lower bounds appearing in (\ref{Dim-U-rank}) and
comparing with Estimate (\ref{CRs-T-U-rank}), we get a more quantitative
form of the general pattern that was hinted before when we looked on Figure %
\ref{cr-vs-dim-gl8}. In particular, for irreps of $U$-rank $k<\left\lfloor 
\frac{n}{2}\right\rfloor $:

\begin{itemize}
\item the dimensions vary by a multiple of $q^{\frac{k(k-1)}{2}}$, although
their CRs are practically the same, of size around $\frac{1}{q^{k}}$;

and,

\item for $n>\frac{(k+1)(k+2)}{2}$, the upper bound for the dimension of $U$%
-rank $k$ irreps is (for sufficiently large $q$) smaller than the lower
bound for rank $k+1$.

But,

\item when $n<\frac{(k+1)(k+2)}{2}$, the range of dimensions for $U$-rank $k$
irreps overlaps (for large enough $q$) the range for $k+1$, and the overlap
grows with $k$. For $k$ in this range, representations of the same dimension
can have different character ratios, which are accounted for by looking at
rank.\smallskip
\end{itemize}

\subsection{\textbf{Cardinality of the Set of }$U$\textbf{-rank }$k$\textbf{%
\ Irreps}}

The cardinality of the set $(\widehat{GL}_{n})_{U,k},$ of irreps of $GL_{n}$
of $U$-rank $k$ can be estimated explicitly. \textbf{%
\begin{figure}[h]\centering
\includegraphics
{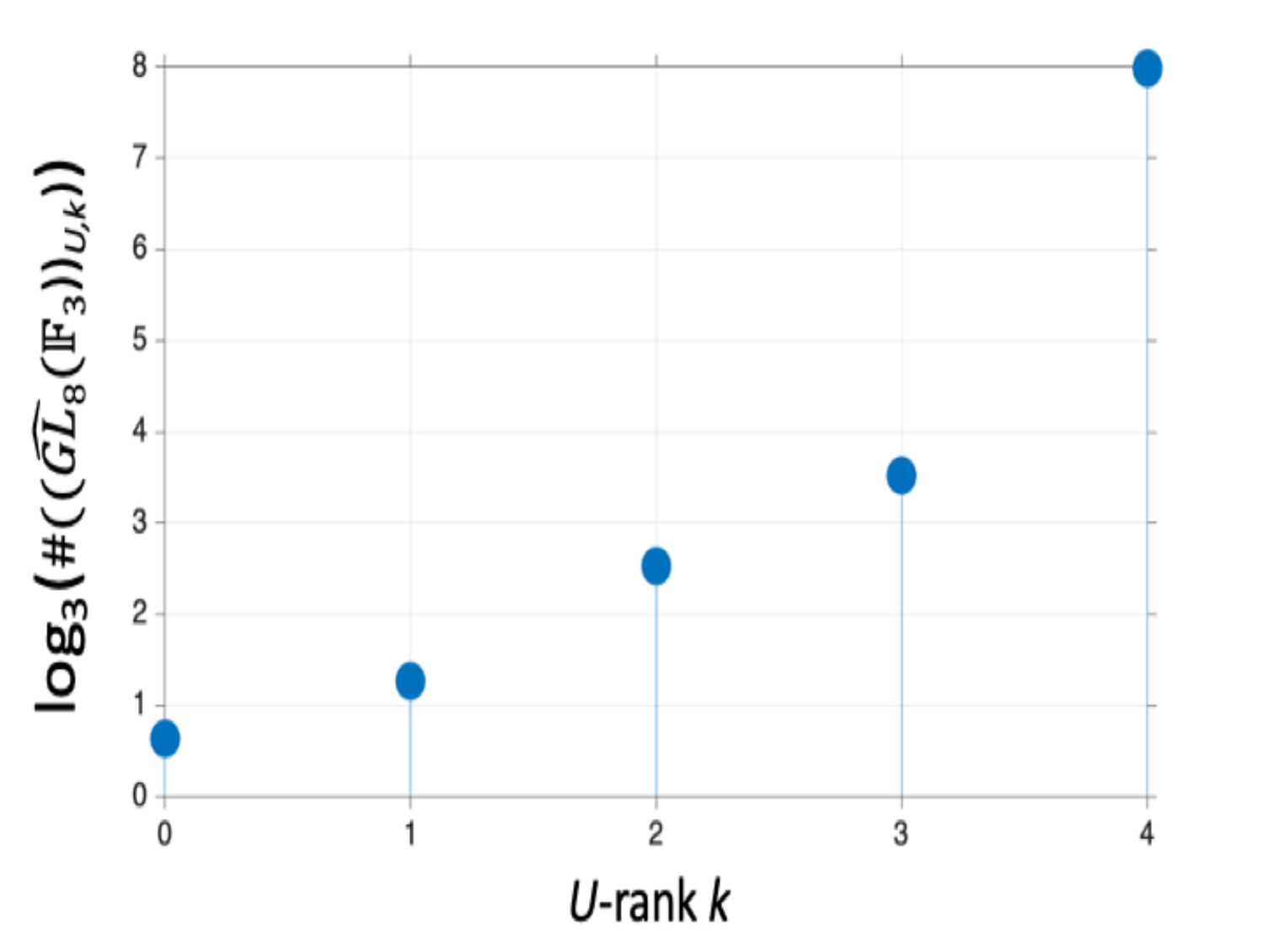}%
\caption{$\log _{q}$-scale of the number of $U$-rank $k$ irreps of $GL_{8}(%
\mathbb{F}_{q}),$ $q=3.$}\label{log-Irr-UR-GL8-3}%
\end{figure}%
}Figure \ref{log-Irr-UR-GL8-3} illustrates how the cardinality of the set $%
(\widehat{GL}_{8}(\mathbb{F}_{3}))_{U,k}$ of $U$-rank $k$ irreps of $GL_{8}(%
\mathbb{F}_{3}),$ grows with $k$.

In general, the following hold:

\begin{theorem}
\label{T-C-URk}We have,%
\begin{equation*}
\#((\widehat{GL}_{n})_{U,k})=\left\{ 
\begin{array}{c}
q^{k+1}+o(...),\text{ \ if }k<\left\lfloor \frac{n}{2}\right\rfloor ; \\ 
q^{n}+o(...),\text{ \ if }k=\left\lfloor \frac{n}{2}\right\rfloor ;%
\end{array}%
\right.
\end{equation*}
\end{theorem}

Our main tool to construct and analyze irreps of $GL_{n}$ of each given $U$%
-rank $k<\left\lfloor \frac{n}{2}\right\rfloor $ is what we discuss next.

\section{\textbf{The eta Correspondence and }$U$\textbf{-Rank\label{S-EC-UR}}%
}

We want to answer the following:\vspace{0.03in}

\textbf{Question: }How to get the information (e.g., the results of Section %
\ref{S-AI-UR}) on $U$-rank $k<\left\lfloor \frac{n}{2}\right\rfloor $ irreps?%
\vspace{0.03in}

In fact, till this point in our story, it is \textit{a priori} not clear why
in general for each $0\leq k\leq \left\lfloor \frac{n}{2}\right\rfloor $
there are at all irreps of $U$-rank $k$?

To answer the above questions for all classical groups we discussed in \cite%
{Gurevich-Howe15, Gurevich-Howe17} the eta correspondence (EC). In the case
of the general linear group $GL_{n},$ it first led in \cite{Gurevich-Howe17}
to an explicit parametrization of the members of the set $(\widehat{GL}%
_{n})_{U,k}$ for $k<\left\lfloor \frac{n}{2}\right\rfloor ,$ enabling the
computation of the cardinality of that set, and in particular to an answer
for the second question above. Secondly, the EC leads to an explicit
Harish-Chandra's "philosophy of cusp forms" type formula (developed in \cite%
{Gurevich-Howe17}, and which will be recalled in Section \ref{S-EC-TR}) for
a general member of $(\widehat{GL}_{n})_{U,k}$ for $k<\left\lfloor \frac{n}{2%
}\right\rfloor .$ This particular formula allowed us in \cite%
{Gurevich-Howe19} to obtain the CRs on the transvection and the dimensions
for the irreps of $U$-rank $k<\left\lfloor \frac{n}{2}\right\rfloor $,
announced in the previous section.

Next we go over some of the details of the basic construction of the EC for $%
GL_{n}$ given in \cite{Gurevich-Howe17}.

\subsection{\textbf{The eta Correspondence\label{S-T-EC}}}

Consider the vector space $L^{2}(M_{n,k})$ of complex valued functions on
the set $M_{n,k}$ of $n\times k$ matrices over $\mathbb{F}_{q}$. The natural
actions of $GL_{n}$ (from the left) and $GL_{k}$ (from the right) on such
matrices induce the pair of commuting actions%
\begin{equation*}
GL_{n}\curvearrowright L^{2}(M_{n,k})\curvearrowleft GL_{k},
\end{equation*}%
that form a single representation%
\begin{equation}
\omega _{n,k}\text{ of }GL_{n}\times GL_{k}\text{ on }L^{2}(M_{n,k}),
\label{omega_n,k}
\end{equation}%
given by $\left[ \omega _{n,k}(g,h)f\right] (m)=f(g^{-1}mh)$, for every $%
h\in GL_{n},$ $m\in M_{n,k},$ $g\in GL_{n}$, and $f\in L^{2}(M_{n,k})$.

The pair $(GL_{n},GL_{k})$ is an example of a dual pair in the language of 
\cite{Howe73}, and $\omega _{n,k}$ is sometime called its \textit{oscillator
(aka Weil) representation}.

Let us decompose $\omega _{n,k}$ into a direct sum of isotypic components
for the irreps of $GL_{k}$,%
\begin{equation}
\omega _{n,k}\simeq \sum_{\tau \in \widehat{GL}_{k}}\overset{\text{rep of }%
GL_{n}\text{ }}{\overbrace{\mathcal{M}(\tau )}}\otimes \tau ,
\label{Isot-Dec}
\end{equation}%
where $\mathcal{\mathcal{M}}(\tau )$ denotes the multiplicity space $%
Hom_{GL_{k}}(\tau ,\omega _{n,k})$ which is a representation of $GL_{n}$.%
\vspace{0.02in}

What can be said about $\mathcal{M}(\tau )$?\vspace{0.02in}

It turns out that for $k\leq \left\lfloor \frac{n}{2}\right\rfloor ,$
although it might be reducible, each $\mathcal{M}(\tau )$ has a unique big
irreducible chunk which is (in a quantitative sense) most of it, and that
can be effectively analyzed. In fact, the notion of $U$-rank gives a way to
distinguish it inside $\mathcal{M}(\tau )$. All of this is the content of
the following:

\begin{theorem}
\label{T-EC-U}Assume $k\leq \left\lfloor \frac{n}{2}\right\rfloor $. We have,

\begin{enumerate}
\item $U$\textbf{-rank }$k$\textbf{\ piece: }Each $\mathcal{\mathcal{M}}%
(\tau )$ contains a unique irreducible component $\eta (\tau )$ of $U$-rank $%
k$, and it appears with multiplicity one, in addition to irreps of lower $U$%
-rank, i.e.,%
\begin{equation*}
\mathcal{M}(\tau )=\overset{U\text{-rank }k}{\overbrace{\eta (\tau )}}\text{
\ }+\text{ \ lower }U\text{-rank irreps.}
\end{equation*}

Moreover,

\item \textbf{eta correspondence}\textit{: The assignment}$\mathit{\ }$\ $%
\tau \mapsto \eta (\tau ),$ defines a one-to-one map \textit{\ }%
\begin{equation}
\eta :\widehat{GL}_{k}\hookrightarrow (\widehat{GL}_{n})_{U,k}.
\label{etaC-U}
\end{equation}%
We call the map (\ref{etaC-U}) the \underline{eta correspondence}.
\end{enumerate}
\end{theorem}

A proof of Parts (1) and (2) of Theorem \ref{T-EC-U} appeared in \cite%
{Gurevich-Howe17}.

We proceed to show that, in the case $k<\left\lfloor \frac{n}{2}%
\right\rfloor ,$ we can say a bit more.

\subsection{\textbf{Exhaustivity of the eta Correspondence\label{S-Ex-EC}}}

Suppose $\rho $ is a representation of $GL_{n}$ and $\chi $ a character
(i.e., one dimensional rep) of this group. We will call the representation $%
\chi \otimes \rho ,$ the \textit{twist of }$\rho $ by $\chi $. Note that
since every character of $GL_{n}$ is trivial on $U$, then the set $(\widehat{%
GL}_{n})_{U,k}$ is preserved under twists by characters.

\begin{example}[All $U$-rank $k=1$ irreps?]
\label{Ex-U-rank-k=1}The oscillator representation $\omega _{n,1}$ of $%
GL_{n}\times GL_{1}$ is given by the natural actions on the space $L^{2}(%
\mathbb{F}_{q}^{n})$ of complex valued functions on the set of column
vectors of length $n$ over $\mathbb{F}_{q}$. Let us assume that $n\geq 4$.
In this case the decomposition (\ref{Isot-Dec}) is%
\begin{equation*}
\omega _{n,1}=\sum_{\lambda \in \widehat{GL}_{1}}\mathcal{M}(\lambda ),
\end{equation*}%
where $\mathcal{M}(\lambda )=\{f:\mathbb{F}_{q}^{n}\rightarrow 
\mathbb{C}
^{\ast }$; $f(av)=\lambda (a)f(v)$, $\ a\in \mathbb{F}_{q}^{\ast },$ $v\in 
\mathbb{F}_{q}^{n}\}.$

It is not difficult to see using direct calculations that,

\begin{enumerate}
\item for $\lambda \neq \mathbf{1}$ the space $\mathcal{M}(\lambda )$ is
irreducible as a $GL_{n}$-representation, it has dimension around $q^{n-1},$
and its CR on $T$ (\ref{T}) is around $\frac{1}{q}$.

In particular, each multiplicative character $\lambda \neq \mathbf{1}$ of $%
GL_{1}=\mathbb{F}_{q}^{\ast }$, is assigned by the EC (\ref{etaC-U}) to the $%
U$-rank $k=1$ irrep 
\begin{equation*}
\eta (\lambda )=\mathcal{M}(\lambda ).
\end{equation*}%
and\vspace{0.01in},

\item The space $\mathcal{M}(\mathbf{1})=(2\times trivial$ $rep)\oplus 
\mathcal{M}(\mathbf{1})_{0}$, where $\mathcal{M}(\mathbf{1})_{0}=\{f\in 
\mathcal{M}(\mathbf{1});$ $f(0)=0$ and $\tsum\limits_{v\in \mathbb{F}%
_{q}^{n}}f(v)=0\}$ is irreducible as a $GL_{n}$-representation, it has
dimension around $q^{n-1},$ and its CR on $T$ is around $\frac{1}{q}$.

In particular, we have the $U$-rank $k=1$ irrep%
\begin{equation*}
\eta (\mathbf{1})=\mathcal{M}(\mathbf{1})_{0}.
\end{equation*}
\end{enumerate}
\end{example}

It can be shown that, if we twist the above irreps by characters we obtain
overall a collection of pairwise non-isomorphic $U$-rank $k=1$ irreps, and
the question is whether we exhausted the set $(\widehat{GL}_{n})_{U,1}$?

In \cite{Gurevich-Howe15, Gurevich-Howe17} we conjectured that the answer to
the above question is yes, and formulated the following:

\begin{conjecture}[Exhaustion]
\label{C-Exhaustion}Suppose $k<\left\lfloor \frac{n}{2}\right\rfloor $.
Then, up to twist by a character, every irrep of $U$-rank $k$ of $GL_{n}$,
is in the image of eta correspondence (\ref{etaC-U}).
\end{conjecture}

In Section \ref{S-P-AC/EC}, we will show that for sufficiently large $q$,
Conjecture \ref{C-Exhaustion} holds true.

\subsection{\textbf{Concluding Remarks on eta and the Analytic Information}}

We would like to remark that,\vspace{0.03in}

\textbf{(A)} \textit{Concerning CRs: }it was shown in \cite{Gurevich-Howe19}
(see also Section \ref{S-EC-TR}) that the description of the irrep $\eta
(\tau )$'s appearing in Part (1) of Theorem \ref{T-C-URk} can be made
effective so one can compute their CRs on the transvection and obtain
Theorem \ref{T-CRs-URank}.

Overall, note for any irrep $\rho \in (\widehat{GL}_{n})_{U,k}$, $%
k<\left\lfloor \frac{n}{2}\right\rfloor ,$ we indeed have%
\begin{equation*}
\frac{\chi _{\rho }(T)}{\dim (\rho )}=\frac{\widehat{1}_{\mathcal{O}_{k}}(T)%
}{\#(\mathcal{O}_{k})}+o(...),
\end{equation*}%
supporting the intuition we had in the process of giving the formal
definition of $U$-rank in Section \ref{S-UR}.\vspace{0.03in}

\textbf{(B) }\textit{Concerning dimensions: }in \cite{Gurevich-Howe19} we
gave an effective description of the $\eta (\tau )$'s, that appears in the
eta correspondence, which in particular, implies, 
\begin{equation*}
\dim (\eta (\tau ))=\dim (\tau )\#(\mathcal{O}_{k})+o(...),
\end{equation*}%
and so, together with the fact (see Lemma \ref{L-M,Mo}) that $\#(\mathcal{O}%
_{k})=q^{k(n-k)}+o(...)$, we see that 
\begin{equation*}
\dim (\eta (\tau ))=\dim (\tau )\cdot q^{k(n-k)}+o(...).
\end{equation*}%
Moreover, the irreps of $GL_{k}$ all have dimensions in the range $1$ to $q^{%
\frac{k(k-1)}{2}}+o(...)$ \cite{Green55}, and we find that Theorem \ref%
{T-Dim-U-rank} follows.\vspace{0.03in} In Section \ref{S-EC-TR}, we recall
another argument from \cite{Gurevich-Howe19} that verifies Theorem \ref%
{T-Dim-U-rank}.\qquad

\textbf{(C) }\textit{Concerning cardinality:} Part (2) of Theorem \ref%
{T-EC-U} combined with the surjectivity of $\eta $ (\ref{etaC-U}) gives the
cardinality of $(\widehat{GL}_{n})_{U,k}$ announced in Theorem \ref{T-C-URk}%
. Indeed, the number of irreps of $GL_{k}$ (the size of the set of conjugacy
classes of that group) is $q^{k}+o(...),$ and in \cite{Gurevich-Howe19} we
showed that after we twist the members in the image of $\eta $ by (the $q-1$%
) characters of $GL_{n}$, you get $q^{k+1}+o(...)$ non-isomorphic irreps, as
claimed.

\section{\textbf{Character Ratios and Tensor Rank}}

Most of the irreps of $GL_{n}$ are of the maximal possible $U$-rank, i.e., $%
\left\lfloor \frac{n}{2}\right\rfloor $. Although the CRs of these irreps
might be relatively small---maybe even too small to contribute to the
harmonic analysis sums, such as (\ref{N-T-Est}), that come up in important
counting problems---it is still the case that we need to estimate them.

To say that an irrep of $GL_{n}$ is of $U$-rank $\left\lfloor \frac{n}{2}%
\right\rfloor $ does an injustice to it from the analytic perspective. For
example, look on the numerical data collected for the group $GL_{8}(\mathbb{F%
}_{3})$ and appear in Figure \ref{cr-vs-ur-gl8}. It shows a large variation
of the CRs at the transvection for the irreps of $U$-rank equal to $%
\left\lfloor \frac{8}{2}\right\rfloor =4$.

So, we want an extension of the notion of $U$-rank in order to control the
CRs $\frac{\chi _{\rho }(T)}{\dim (\rho )}$ also within the irreps of $U$%
-rank $k=\left\lfloor \frac{n}{2}\right\rfloor $.

In \cite{Gurevich-Howe17} we proposed such an extension, called tensor rank,
for representations of all classical groups. In fact it appeared with
different terminology already in the unpublished notes \cite{Howe73}. We
proceed to discuss this notion in the case of $GL_{n}$, where we showed in 
\cite{Gurevich-Howe19} that it does a pretty good job---see Figure \ref%
{cr-vs-tr-gl8}.

\subsection{\textbf{Tensor Rank: Formal Definition and Agreement with }$U$%
-Rank\textbf{\label{S-TR}}}

The definition of tensor rank will be given using an extension of the way we
realized the set $(\widehat{GL}_{n})_{U,k},$ for $k<\left\lfloor \frac{n}{2}%
\right\rfloor $.

Consider the oscillator representation\footnote{%
Up to a sign, $\omega _{n}$ is the restriction of the oscillator
representation of $Sp_{2n}$ to $GL_{n}$ \cite{Gerardin77, Howe73, Weil64}.} $%
\omega _{n}$ of $GL_{n}$ given by its right action on the space of complex
valued functions $L^{2}(\mathbb{F}_{q}^{n})$ on $\mathbb{F}_{q}^{n}$, or
more generally consider its $k$-fold tensor product%
\begin{equation*}
GL_{n}\overset{\omega _{n}^{\otimes ^{k}}}{\curvearrowright }L^{2}(M_{n,k})%
\text{,}
\end{equation*}%
given using the right action of this group on $n\times k$ matrices. Note
that $\omega _{n}^{\otimes ^{k}}$is just the restriction of $\omega _{n,k}$ (%
\ref{omega_n,k}) to $GL_{n}$.

Denote by $\widehat{GL}_{n}(\omega _{n}^{\otimes ^{k}})$ the set of irreps
of $GL_{n}$ that appear in $\omega _{n}^{\otimes ^{k}},$ and by $\mathbf{1}$
the trivial representation. In \cite{Gurevich-Howe19} we showed that,

\begin{proposition}
We have a sequence of proper containments%
\begin{equation}
\{\mathbf{1}\}\subsetneqq \widehat{GL}_{n}(\omega _{n}^{\otimes
^{1}})\subsetneqq \ldots \subsetneqq \widehat{GL}_{n}(\omega _{n}^{\otimes
^{n}})=\widehat{GL}_{n}.  \label{TRF}
\end{equation}
\end{proposition}

Now, looking at (\ref{TRF}) and taking into account the action of characters
(i.e., $1$-dim representations) on irreps, we introduce in \cite%
{Gurevich-Howe17} the following:

\begin{definition}[\textbf{Tensor rank}]
\label{D-TR}We will say that $\rho \in \widehat{GL}_{n}$ is of \underline{%
\textbf{tensor rank}} $k$, if the minimal $\ell $ that we can write it as a
tensor product of a character and an irrep from $\widehat{GL}_{n}(\omega
_{n}^{\otimes ^{\ell }})$ is $\ell =k$.
\end{definition}

We may use the notations $\otimes $-$rank(\rho )=k$, or $rank_{\otimes
}(\rho )=k$, to indicate that a representation $\rho $ of $GL_{n}$ has
tensor rank $k$, and denote the set of all such irreps by $(\widehat{GL}%
_{n})_{\otimes ,k}$.

\begin{remark}
The notion of tensor rank induces (and is defined by) a filtration on the
representation ring 
\begin{equation*}
R(GL_{n})=%
\mathbb{Z}
\lbrack \widehat{GL}_{n}],
\end{equation*}%
generated from the set $\widehat{GL}_{n}$ using the operations of addition
and multiplication given, respectively, by direct sum $\oplus $ and tensor
product $\otimes .$ Indeed, let us extend the definition of tensor rank to
arbitrary (not necessarily irreducible) representation of $GL_{n}$ and say
it is of tensor rank $k$ if it contains irreps of tensor rank $k$ but not of
higher tensor rank. In particular, we have the \underline{\textit{tensor
rank filtration}} which is obtained by letting $F_{\otimes ,k}$ be the
collection of elements of $R(G)$ that are sums of irreps of tensor rank less
or equal to $k,$ and it satisfies

\begin{itemize}
\item $F_{\otimes ,(k-1)}\subset F_{\otimes ,k}$, $F_{\otimes ,i}\otimes
F_{\otimes ,j}\subset F_{\otimes ,i+j}$ for every $i,j,k;$

and

\item $F_{\otimes ,n}=R(G)$.
\end{itemize}
\end{remark}

\vspace{0.01in}

Finally, note that Part (1) of Theorem \ref{T-EC-U} implies that 
\begin{equation}
(\widehat{GL}_{n})_{\otimes ,k}\subset (\widehat{GL}_{n})_{U,k},\text{ \ for
\ }k<\left\lfloor \frac{n}{2}\right\rfloor ,  \label{TR-UR}
\end{equation}%
and Conjecture \ref{C-Exhaustion} can be restated as follows:

\begin{conjecture}[Agreement]
\label{Con-AC}The inclusion (\ref{TR-UR}) should be replaced by equality.
\end{conjecture}

In particular, we conclude that, tensor rank is a natural extension of the
notion of $U$-rank. But, is it going to split nicely the collection of $U$%
-rank $k=\left\lfloor \frac{n}{2}\right\rfloor $ irreps of $GL_{n}$?

\subsection{\textbf{Numerics for Character Ratios vs. Tensor }Rank\textbf{\ 
\label{S-N-CRvsTRank}}}

\begin{figure}[h]\centering
\includegraphics
{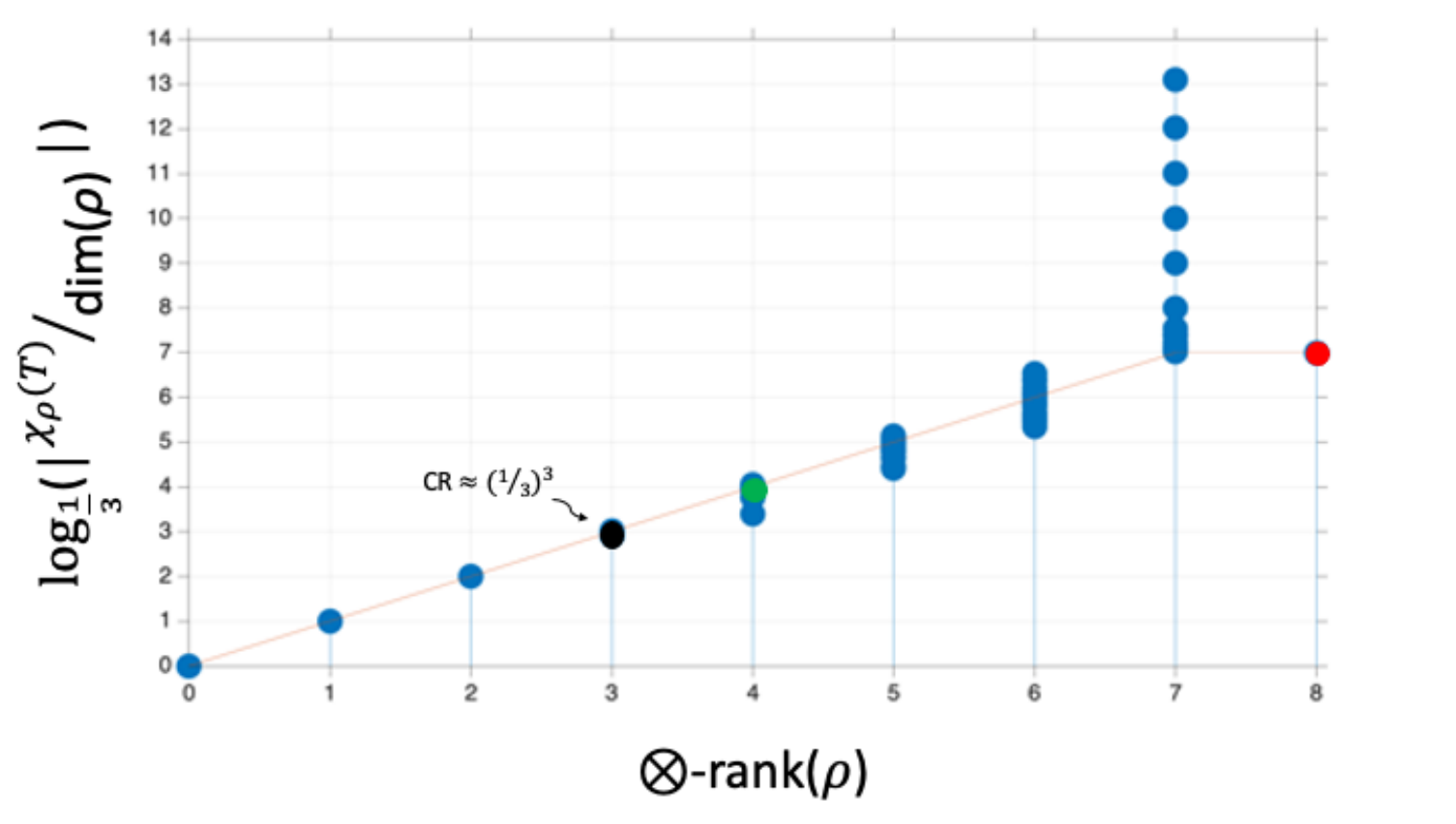}%
\caption{$\log _{\frac{1}{q}}$-scale of CRs\ vs. $\otimes $-rank for irreps $%
\protect\rho $ of $GL_{8}(\mathbb{F}_{q}),$ $q=3.$}\label{cr-vs-tr-gl8}%
\end{figure}

The answer to the above question seems to be yes and, before we write down
formal statements, we would like to see this with the aid of some supportive
numerical data collected for the group $GL_{8}(\mathbb{F}_{3})$ which
appears in Figures \ref{cr-vs-tr-gl8}, \ref{cr-vs-ur-gl8}, and \ref%
{cr-vs-dim-gl8}. 

Recall that Figure \ref{cr-vs-ur-gl8} illustrates the general fact that for
irreps of $GL_{n}$ of $U$-rank $k<\left\lfloor \frac{n}{2}\right\rfloor $,
the CRs at the transvection $T$ are all essentially of the same size $\frac{1%
}{q^{k}}$, despite the fact that the dimensions of the irreps involved might
vary by a multiple of large power of $q$. Due to the agreement conjecture
this should also holds---and illustrated in Figure \ref{cr-vs-tr-gl8}---for
low $\otimes $-rank irreps, i.e., these of $\otimes $-rank $k<\left\lfloor 
\frac{n}{2}\right\rfloor $.

The next thing that Figure \ref{cr-vs-tr-gl8} illustrates is that indeed
(compare with Figure \ref{cr-vs-ur-gl8}) the tensor rank invariant splits
further the collection of irreps of $U$-rank of $k=\left\lfloor \frac{n}{2}%
\right\rfloor $, and that the splitting seems to do more or less the job we
wanted the tensor rank to do. Specifically, Figure \ref{cr-vs-tr-gl8}
illustrates the fact that for tensor rank $\frac{n}{2}\leq k\leq n-1,$ the
CRs at the transvection are of the order of magnitude of $\frac{1}{q^{k}}$
time a constant (independent of $q$), and that for all tensor rank $n$
irreps the CRs are exactly $\frac{1}{q^{n-1}}$ in absolute value. Finally, a
look at the black-green-red circles above 16 in Figure \ref{cr-vs-dim-gl8},
and how they appear in Figure \ref{cr-vs-tr-gl8}, illustrates the fact that
tensor rank provides a reason for why irreps of the same order of magnitude
of dimension can have very different CRs at $T$, namely, the answer is that
these irreps have different tensor ranks.

The above numerical results can be quantified precisely and proved. This is
part of what we do next.

\section{\textbf{Analytic Information on Tensor Rank }$k$\textbf{\ Irreps of 
}$GL_{n}$\label{S-AI-TRk-GLn}}

In this section we present results obtained in \cite{Gurevich-Howe19}
concerning the character ratios and dimensions of the irreps of $\otimes $%
-rank $k,$ i.e., the members of $(\widehat{GL}_{n})_{\otimes ,k},$ including
the cardinality of that set. Although we will not prove these results in
this note, for the benefit of the reader we explain in Section \ref{S-EC-TR}
what are the main sources of informations that enable us to derive them.

\subsection{\textbf{Character Ratios on the Transvection\label{S-CR-T-TR}}}

For the CRs on the transvection $T$ (\ref{T}), the following, essentially
sharp, estimate in term of the tensor rank holds:

\begin{theorem}
\label{T-AI-GLn}Fix $0\leq k\leq n$. Then, for $\rho \in (\widehat{GL}%
_{n})_{\otimes ,k},$ we have an estimate:

\begin{equation}
\frac{\chi _{\rho }(T)}{\dim (\rho )}=\left\{ 
\begin{array}{c}
\text{ }\frac{1}{q^{k}}+o(...)\text{, \ \ \ \ if \ \ \ }k<\frac{n}{2};\text{
\ \ \ \ \ \ \ \ \ \ \ \ } \\ 
\text{\ \ \ \ \ \ } \\ 
\frac{c_{\rho }}{q^{k}}+o(...)\text{, \ \ \ \ if \ \ }\frac{n}{2}\leq k\leq
n-1;\text{\ \ } \\ 
\text{\ } \\ 
\text{ \ }\frac{-1}{q^{n-1}-1}\text{, \ \ \ \ \ \ \ \ \ \ if \ \ \ }k=n,%
\text{ \ \ \ \ \ \ \ \ \ \ }%
\end{array}%
\right.  \label{CRs-GLn}
\end{equation}%
where $c_{\rho }$ is a certain integer (independent of $q$) combinatorially
associated with $\rho $.
\end{theorem}

\begin{remark}
For irreps $\rho $ of tensor rank $\frac{n}{2}\leq k\leq n-1,$ the constant $%
c_{\rho }$ in (\ref{CRs-GLn}) might be equal to zero. In this case, the
estimate on $\frac{\chi _{\rho }(T)}{\dim (\rho )}$ is simply $o(\frac{1}{%
q^{k}}).$ However, the possibility of $c_{\rho }=0$ is fairly rare, and (at
least for $k\neq n-1$) we are not sure if it happens at all.
\end{remark}

Note that (\ref{CRs-GLn}) is a formal validation to some of the phenomena
that Figure \ref{cr-vs-tr-gl8} illustrates.

\subsection{\textbf{Dimensions of Irreps\label{S-Dim-Trank}}}

\begin{figure}[h]\centering
\includegraphics
{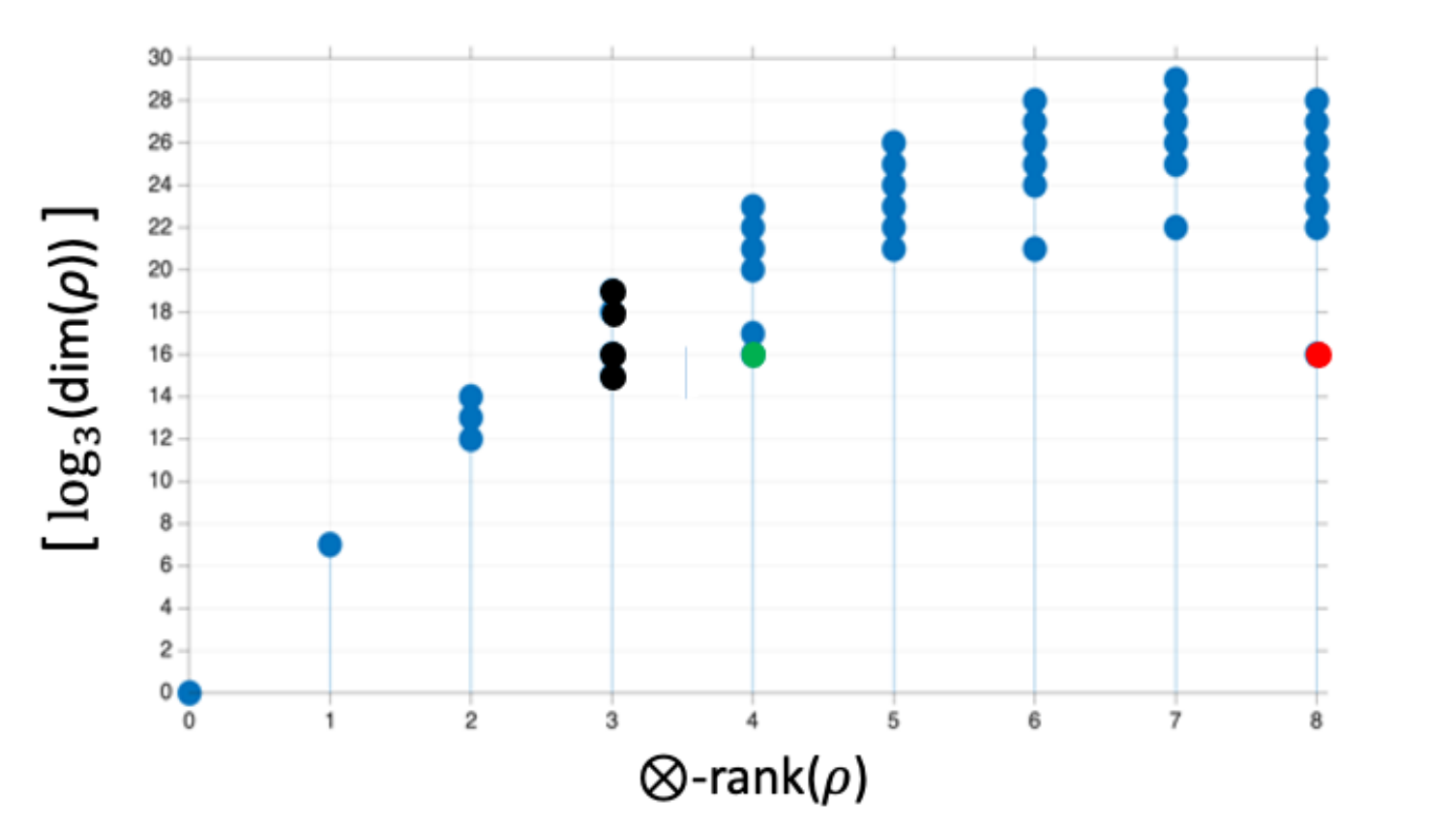}%
\caption{$\log _{q}$-scale of dimension vs. $\otimes $-rank for irreps $%
\protect\rho $ of $GL_{8}(\mathbb{F}_{q}),$ $q=3.$}\label{dim-vs-tr-gl8}%
\end{figure}

We proceed to present information on the dimensions of the irreps of tensor
rank $k$. Figure \ref{dim-vs-tr-gl8} gives a numerical illustration for the
distribution of the dimensions of the irreps of $GL_{8}(\mathbb{F}_{3})$
within each given tensor rank.

The following are the sharp lower and upper bounds obtained in \cite%
{Gurevich-Howe19} (and that formally explain Figure \ref{dim-vs-tr-gl8}; the
black-green-red dots were discussed in Section \ref{S-N-CRvsD}) on the
dimensions of the $\otimes $-rank $k$ irreps:

\begin{theorem}
\label{T-dim-GLn}Fix $0\leq k\leq n$. Then, for $\rho \in (\widehat{GL}%
_{n})_{\otimes ,k},$ we have an estimate:%
\begin{equation}
q^{k(n-k)+\frac{k(k-1)}{2}}+o(...)\geq \dim (\rho )\geq \left\{ 
\begin{array}{c}
q^{k(n-k)}+o(\ldots ),\text{ \ \ \ \ \ \ \ \ \ if\ \ }k<\frac{n}{2};\text{ \
\ \ \ \ \ \ \ \ \ } \\ 
\\ 
\text{ }q^{(n-k)(3k-n)}+o(\ldots ),\text{ \ \ \ if \ }\frac{n}{2}\leq k<%
\frac{2n}{3};\text{ \ \ \ } \\ 
\\ 
\text{ \ \ }q^{k(n-k)+\frac{k^{2}}{4}}+o(\ldots ),\text{ \ \ \ \ \ if \ }%
\frac{2n}{3}\leq k\leq n,\text{ even}; \\ 
\\ 
\text{\ \ \ \ \ }q^{k(n-k)+\frac{(k-3)^{2}}{4}+3(k-2)}+o(...),\text{ \ if }%
\frac{2n}{3}\leq k\leq n,\text{ odd;\ \ }%
\end{array}%
\right.  \label{Dim-GLn}
\end{equation}%
Moreover, the upper and lower bounds in (\ref{Dim-GLn}) are attained.
\end{theorem}

\subsection{\textbf{The Number of Irreps of Tensor Rank }$k$\textbf{\ of }$%
GL_{n}$\label{S-Card-Trank-k}}

Finally, we present information concerning the cardinality of the set of
irreps of $\otimes $-rank $k$---see Figure \ref{log-irr(g)_tr,k-gl8} for
illustration.

\textbf{%
\begin{figure}[h]\centering
\includegraphics
{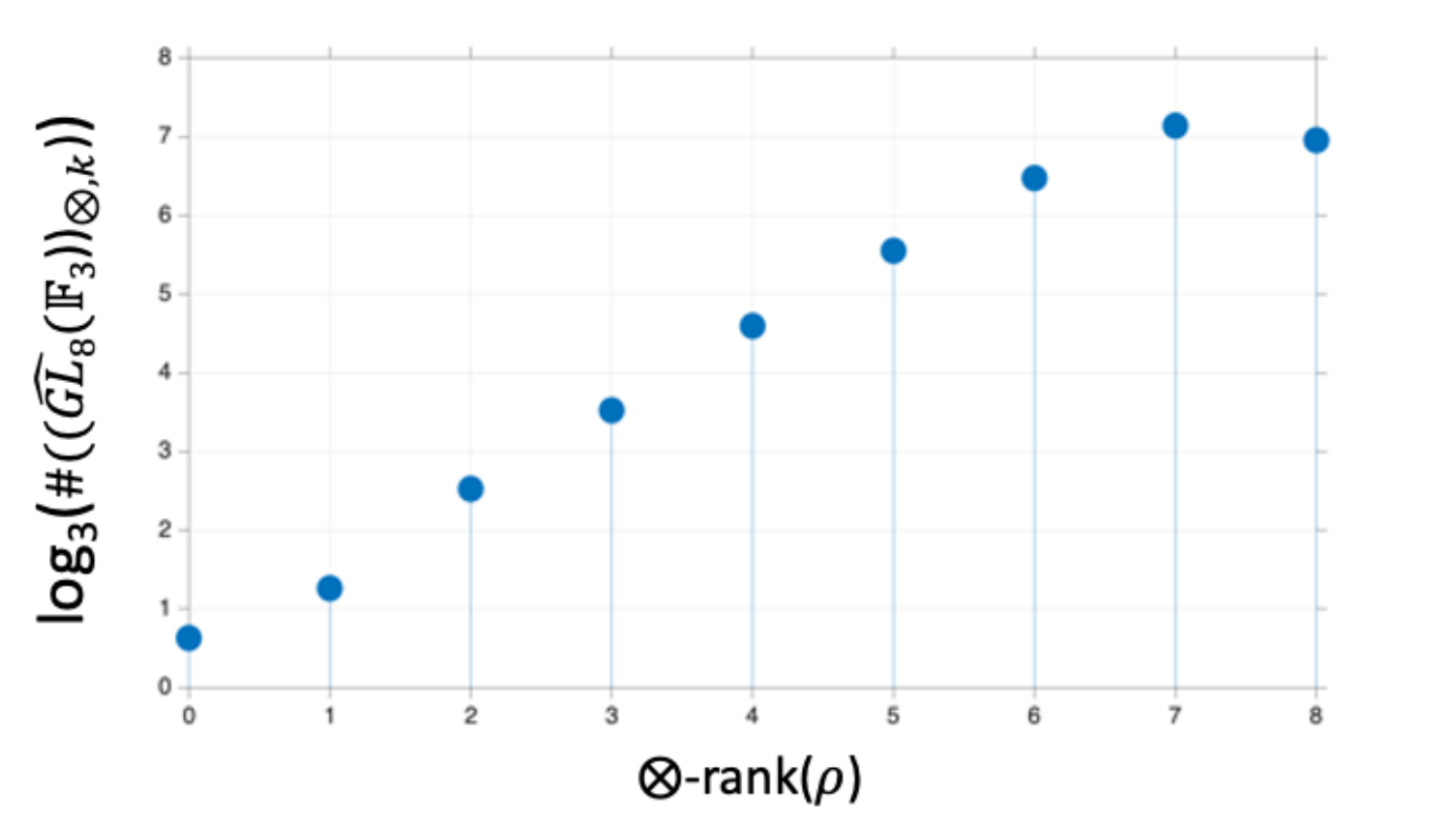}%
\caption{$\log _{q}$-scale of the number of $\otimes $-rank $k$ irreps of $%
GL_{8}(\mathbb{F}_{q}),$ $q=3.$}\label{log-irr(g)_tr,k-gl8}%
\end{figure}%
}

In this aspect, we have the following essentially sharp estimate:

\begin{theorem}
\label{T-Card-trank-k}Fix $0\leq k\leq n$. Then, we have, 
\begin{equation}
\#((\widehat{GL}_{n})_{\otimes ,k})=\left\{ 
\begin{array}{c}
q^{k+1}+o(...)\text{, \ if \ }k\leq n-2; \\ 
c_{k}q^{n}+o(...)\text{, \ if \ }n-2<k,%
\end{array}%
\right.  \label{Card-k-GLn}
\end{equation}%
where $0<c_{n-1},$ $c_{n}<1,$ $c_{n-1}+c_{n}=1$.
\end{theorem}

\subsection{\textbf{Some Remarks}}

We would like to make several remarks concerning the analytic information
given just above, that extend in a bit more detailed way similar remarks
given for $U$-rank in Section \ref{S-Dim-Urank-k}.

\subsubsection{\textbf{Tensor Rank vs. Dimension as Indicator for Size of
Character Ratio\label{S-TRvsDim-for-CR}}}

Looking back on the analytic information presented in the Sections \ref%
{S-CR-T-TR} and \ref{S-Dim-Trank}, we observe the following:\smallskip

\textbf{(A) For irreps in a given tensor rank.\smallskip }

\begin{figure}[h]\centering
\includegraphics
{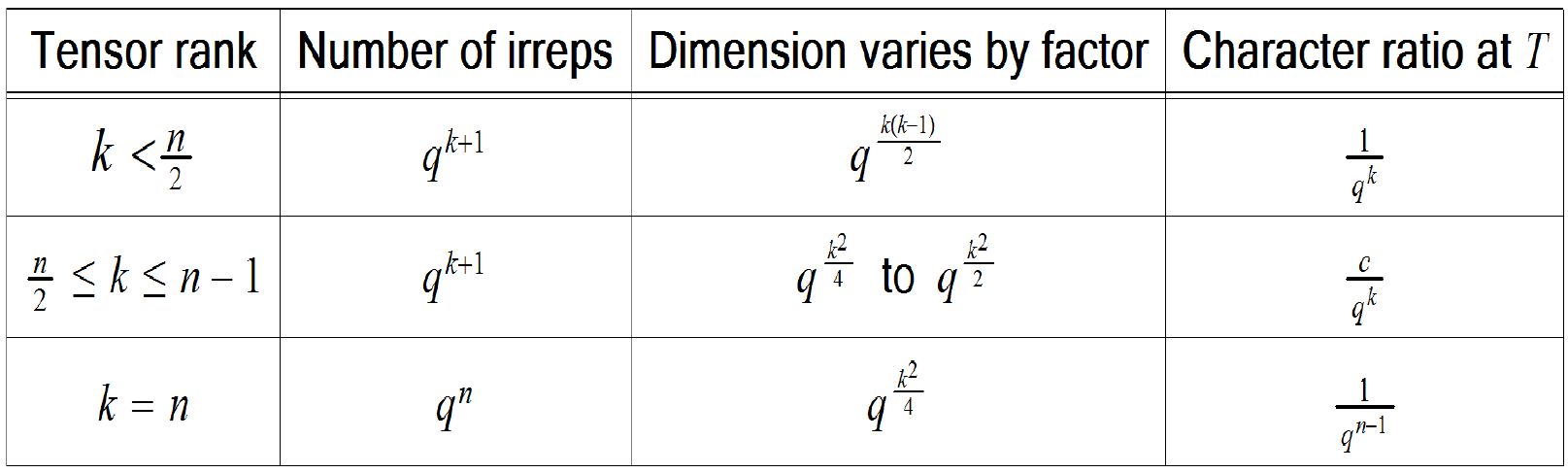}%
\caption{CRs vs. variation in dimensions (in order of magnitude of power of $%
q$) for $\protect\widehat{GL}_{n}$. }\label{cr-vs-dim-trank-k-gl7_3}%
\end{figure}

A comparison of (\ref{Dim-GLn}) and (\ref{CRs-GLn}) demonstrates---see
Figure \ref{cr-vs-dim-trank-k-gl7_3} for a summary---what we illustrated in
Sections \ref{S-N-CRvsD} and \ref{S-N-CRvsTRank}: Within a given tensor rank 
$k$ the 
dimensions may vary by a large factor (around $q^{\frac{k(k-1)}{2}}$ for
rank $k<\frac{n}{2}$, and between $q^{\frac{k^{2}}{4}}$ to $q^{\frac{k^{2}}{2%
}}$ for $\frac{n}{2}\leq k$ - quantities are given in approximate order of
magnitude of power of $q$) but the CRs at the transvection are practically
the same, of size around $\frac{1}{q^{k}}$ (for $\frac{n}{2}\leq k\leq n-1$
a multiple of $\frac{1}{q^{k}}$ by a constant independent of $q$).

\textbf{(B) For irreps of different tensor ranks.\smallskip }

Looking on (\ref{Dim-GLn}) we notice that:

\begin{itemize}
\item for $n>\frac{(k+1)(k+2)}{2}$, the upper bound for the dimension of $%
\otimes $-rank $k$ irreps is (for sufficiently large $q$) smaller than the
lower bound for rank $k+1$.

But,

\item when $n<\frac{(k+1)(k+2)}{2}$, the range of dimensions for $\otimes $%
-rank $k$ irreps overlaps (for large enough $q$) the range for $k+1$, and
the overlap grows with $k$. For $k$ in this range, representations of the
same dimension can have different character ratios, which are accounted for
by looking at rank.\smallskip
\end{itemize}

In conclusion, it seems that tensor rank of a representation is a better
indicator than dimension for the size of its character ratio, at least on
elements such as the transvection.

\subsubsection{\textbf{Comparison with Existing Formulations in the
Literature}}

In most of the literature on character ratios that we are aware of (see,
e.g., \cite{Bezrukavnikov-Liebec-Shalev-Tiep18} and the papers cited
there.), estimates on character ratios are given in terms of the dimension
of representations.

Although the dimension is a standard invariant of representations, as we
have seen in Parts (A) and (B) of Section \ref{S-TRvsDim-for-CR}, the
dimensions of representations with a given tensor rank can vary
substantially (i.e., by large powers of $q$), while the character ratio at
the transvection stays more or less constant (at least for $k<\frac{n}{2}$).
Thus, using only dimension to bound character ratio might lead to
non-optimal estimates.

\begin{figure}[h]\centering
\includegraphics
{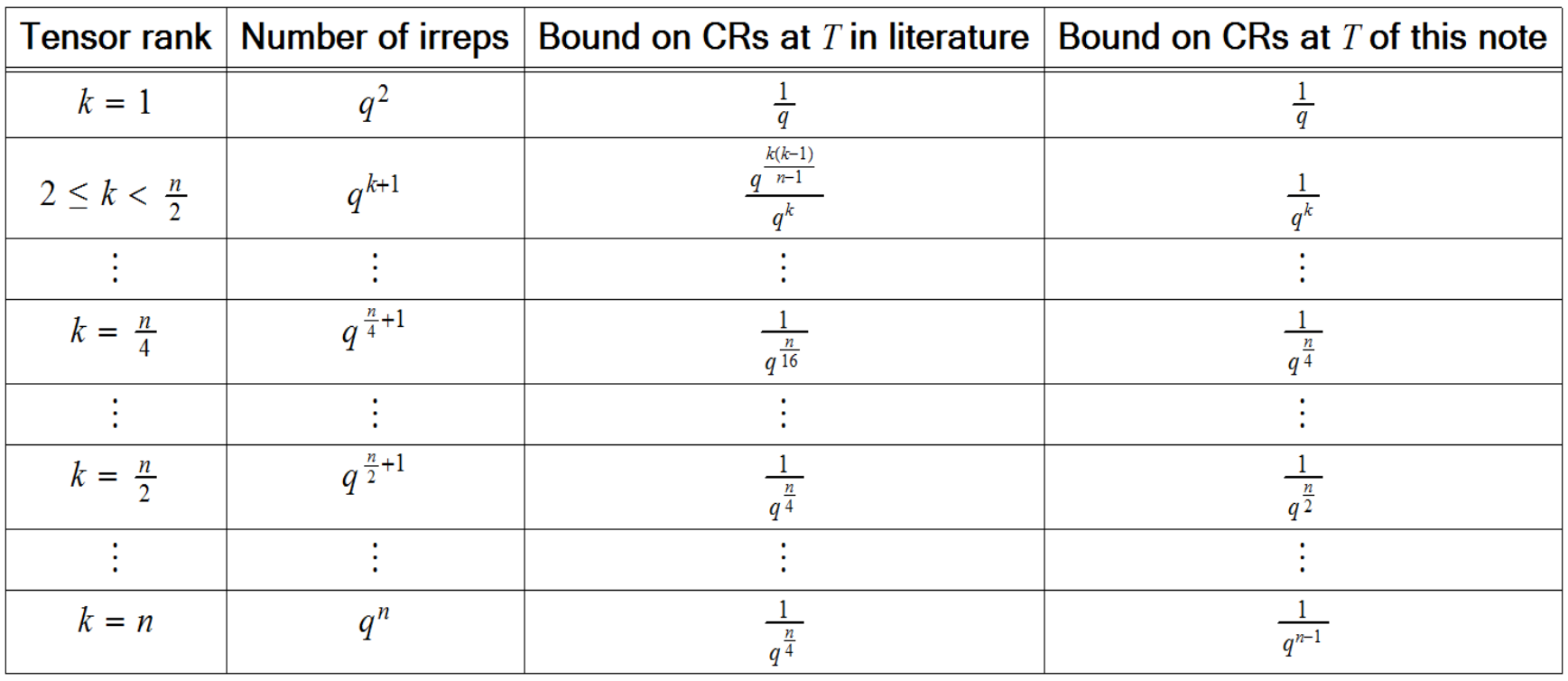}%
\caption{Bounds on CRs: Current literature vs. this note (in order of
magnitude). }\label{compcrs-gln}%
\end{figure}

In particular, the estimates in this note for the character ratio on the
transvection are optimal (in term of the tensor rank), and are, in general,
stronger than the corresponding estimates in the paper cited above. For
example, for $k<\frac{n}{2}$, rather than the bound of $\frac{1}{q^{k}}$,
the paper \cite{Bezrukavnikov-Liebec-Shalev-Tiep18} gives bounds of the
order of magnitude of $\frac{q^{\frac{k(k-1)}{n-1}}}{q^{k}},$ and the
exponent $\frac{k(k-1)}{n-1}$ can be fairly large when $n$ is large and $k$
is near $\frac{n}{2}$. The table in Figure \ref{compcrs-gln} gives some
examples of the relationship between the results of this note, and of the
literature cited above.%

Recall that for the motivational example described in the introduction (see
Section \ref{S-HAofGP}) we wanted to have information on irreps of $%
SL_{n}=SL_{n}(\mathbb{F}_{q}),$ $n\geq 3$. These can be deduced from the one
we just formulated above for $GL_{n}$, as $SL_{n}$ is a very big subgroup of 
$GL_{n}$.

\section{\textbf{Analytic Information on Tensor Rank }$k$\textbf{\ Irreps of 
}$SL_{n}$\textbf{\label{S-AI-SLn}}}

In this section we formulate the analytic results obtained in \cite%
{Gurevich-Howe19} for the tensor rank $k$ irreps of $SL_{n},$ $n\geq 3$. The
case of $SL_{2}$ is somewhat special and \cite{Gurevich-Howe18} was devoted
to its description.

\subsection{\textbf{Tensor Rank for Representations of }$SL_{n}$}

First we introduce the following terminology. We assume $n\geq 3$.

\begin{definition}
\label{D-TR-SLn}We will say that an irreducible representation $\pi $ of $%
SL_{n}$ has \underline{tensor rank} $k$ if it appears in the restriction of
a tensor rank $k$ (and not less) irrep of $GL_{n}.$
\end{definition}

As before, we denote by $(\widehat{SL}_{n})_{\otimes ,k}$ the set of irreps
of $SL_{n}$ of $\otimes $-rank $k$.\smallskip

\begin{remark}
\label{R-TR-STR-SLn}Note that the condition that $\pi $ should satisfy in
Definition \ref{D-TR-SLn} is equivalent to the requirement that (replacing $%
GL_{n}$ by $SL_{n}$) in (\ref{TRF}) it will appear in $\widehat{SL}%
_{n}(\omega _{n}^{\otimes \ell })$ for $\ell =k,$ but not for $\ell <k$.
\end{remark}

The main technique used in \cite{Gurevich-Howe19} to get information on
irreps of $SL_{n}$ is through the way they appear inside irreps of $GL_{n}$.
This can be understood using the Clifford-Mackey's theory \cite{Clifford37,
Mackey49}. We will not repeat the analysis in this note, but the "intuitive
picture" is that nearly every irrep of $GL_{n}$ stays irreducible after
restriction to $SL_{n}$, and hence---using some additional favorite
facts---the results for $SL_{n}$ are the "same" as for $GL_{n}$.

We start with the estimates on the character ratios.

\subsection{\textbf{Character Ratios on the Transvection}}

We have the following sharp estimates in term of the tensor rank:

\begin{corollary}
\label{C-CRs-SLn}Fix $n\geq 3$, and $0\leq k\leq n$. Then, for $\pi \in (%
\widehat{SL}_{n})_{\otimes ,k},$ we have,%
\begin{equation}
\frac{\chi _{\pi }(T)}{\dim (\pi )}=\left\{ 
\begin{array}{c}
\text{ }\frac{1}{q^{k}}+o(...)\text{, \ \ \ \ if \ \ \ }k<\frac{n}{2};\text{
\ \ \ \ \ \ \ \ \ \ \ \ } \\ 
\text{\ \ \ \ \ \ } \\ 
\frac{c_{\pi }}{q^{k}}+o(...)\text{, \ \ \ \ if \ \ }\frac{n}{2}\leq k\leq
n-1;\text{\ \ } \\ 
\text{\ } \\ 
\text{ \ }\frac{-1}{q^{n-1}-1}\text{, \ \ \ \ \ \ \ \ \ \ if \ \ \ }k=n,%
\text{ \ \ \ \ \ \ \ \ \ \ }%
\end{array}%
\right.  \label{CRs-SLn}
\end{equation}%
where $c_{\pi }$ is a certain integer (independent of $q$) combinatorially
associated with $\pi $.
\end{corollary}

\begin{remark}
For irreps $\pi $ of tensor rank $\frac{n}{2}\leq k\leq n-1,$ the constant $%
c_{\pi }$ in (\ref{CRs-SLn}) might be equal to zero. In this case, the
estimate on $\frac{\chi _{\pi }(T)}{\dim (\pi )}$ is simply $o(\frac{1}{q^{k}%
}).$
\end{remark}

\subsection{\textbf{Dimensions of Irreps}}

It turns out that most of the irreps of $GL_{n}$ that give the lower and
upper bounds on dimensions of tensor rank $k$ irreps, stays irreducible as
representations of $SL_{n}$. As a consequence, from the corresponding
results for $GL_{n}$, we obtain,

\begin{corollary}
\label{C-Dim-SLn}Fix $n\geq 3$, and $0\leq k\leq n$. Then, for $\pi \in (%
\widehat{SL}_{n})_{\otimes ,k}$, we have an estimate:%
\begin{equation}
q^{k(n-k)+\frac{k(k-1)}{2}}+o(...)\geq \dim (\pi )\geq \left\{ 
\begin{array}{c}
q^{k(n-k)}+o(\ldots ),\text{ \ \ \ \ \ \ \ \ \ if\ \ }k<\frac{n}{2};\text{ \
\ \ \ \ \ \ \ \ \ } \\ 
\\ 
\text{ }q^{(n-k)(3k-n)}+o(\ldots ),\text{ \ \ \ if \ }\frac{n}{2}\leq k<%
\frac{2n}{3};\text{ \ \ \ } \\ 
\\ 
\text{ \ \ }q^{k(n-k)+\frac{k^{2}}{4}}+o(\ldots ),\text{ \ \ \ \ \ if \ }%
\frac{2n}{3}\leq k\leq n,\text{ even}; \\ 
\\ 
\text{\ \ \ \ \ }q^{k(n-k)+\frac{(k-3)^{2}}{4}+3(k-2)}+o(...),\text{ \ if }%
\frac{2n}{3}\leq k\leq n,\text{ odd;\ \ }%
\end{array}%
\right.  \label{Dim-SLn}
\end{equation}%
Moreover, the upper and lower bounds in (\ref{Dim-SLn}) are attained.
\end{corollary}

\subsection{\textbf{The Number of Irreps of Tensor Rank }$k$\textbf{\ of }$%
SL_{n}$}

The fact, mentioned earlier, that most tensor rank $k$ irreps of $GL_{n}$
stay irreducible after restricting them to $SL_{n}$, is the core fact used
in \cite{Gurevich-Howe19} to deduce (see estimates (\ref{Card-k-GLn})) the
following:

\begin{proposition}
\label{P-Num-Irr-k-SLn}Fix $n\geq 3$, and $0\leq k\leq n$. Then, we have,%
\begin{equation}
\#((\widehat{SL}_{n})_{\otimes ,k})=\left\{ 
\begin{array}{c}
\text{ \ }q^{k}+o(...)\text{, \ \ \ if \ }k\leq n-2; \\ 
c_{k}q^{k-1}+o(...)\text{, \ if \ }n-2<k,%
\end{array}%
\right.  \label{Card-SLn}
\end{equation}%
where $0<c_{n-1},$ $c_{n}<1,$ with $c_{n-1}+c_{n}=1$.
\end{proposition}

\section{\textbf{Back to the Generation Problem\label{S-BtoGP}}}

Having the analytic information on the irreps of $SL_{n}$, $n\geq 3$, we can
address the generation problem discussed in the introduction. In particular,
we can derive Theorem \ref{T-STS}.

\subsection{\textbf{Setting and Statement}}

We considered (see the introduction) the conjugacy class $C\subset SL_{n}$
of the transvection $T$ (\ref{T}), and for $\ell \geq n,$ we looked at the
set%
\begin{equation}
M_{\ell ,g}=\left\{ (c_{1},\ldots ,c_{\ell })\in C^{\ell };\text{ \ }%
c_{1}\cdot \ldots \cdot c_{\ell }=g\right\} ,  \label{Mlg-again}
\end{equation}%
where $g$ is an element of $SL_{n}$ which one can't form by less than $n$
products from $C,$ and is regular semi-simple, i.e., all its eigenvalues
over an algebraic closure of $\mathbb{F}_{q}$ are different. We explained
that this means that our $g$ is semi-simple regular and $1$ is not one of
its eigenvalues. Let us denote the set of all such $g$'s by $(\partial
(G))_{reg}$.

We wanted to compute the cardinality of the set $M_{\ell ,g}$ (\ref%
{Mlg-again}), and to show:\vspace{0.02in}

\textbf{Theorem \ref{T-STS} (Set-theoretic size - restated). }\textit{For an
element }$g\in (\partial (G))_{reg}$, we have,%
\begin{equation}
\#(M_{\ell ,g})=\frac{\#(C^{\ell })}{\#(SL_{n})}\cdot \left\{ 
\begin{array}{c}
1-O(\frac{1}{q}),\text{ \ \ \ \ \ \ \ \ \ \ \ \ \ \ \ \ if \ }\ell =n; \\ 
1-(\frac{2}{q}(\frac{1}{q})^{\ell -n}+o(...)),\text{ \ if \ }\ell >n.%
\end{array}%
\right.  \label{Card-Mlg-again}
\end{equation}

\subsection{\textbf{A Proof of the Set-Theoretic Size Theorem }}

To proof Theorem \ref{T-STS}, we follow the strategy proposed in Section \ref%
{S-HAofGP} and invoke harmonic analysis for our purpose. We have the
Frobenius type formula%
\begin{equation}
\#(M_{\ell ,g})=\frac{\#(C^{\ell })}{\#(SL_{n})}\cdot (1+\overset{S_{\ell ,g}%
}{\overbrace{\underset{\mathbf{1}\neq \pi \in \widehat{SL}_{n}}{\sum }\dim
(\pi )\left( \frac{\chi _{\pi }(T)}{\dim (\pi )}\right) ^{\ell }\chi _{\pi
}(g^{-1})})},  \label{F-C-Mlg-again}
\end{equation}%
and so we just need to show the sum $S_{\ell ,g}$ in (\ref{F-C-Mlg-again})
is of the size of the error term in statement (\ref{Card-Mlg-again}).

At this point we already know about tensor rank and can further split $%
S_{\ell ,g}$ over the various ranks%
\begin{equation}
S_{\ell ,g}=\sum_{k=1}^{n}\overset{(S_{\ell ,g})_{k}}{\overbrace{\underset{%
\pi \in (\widehat{SL}_{n})_{\otimes ,k}}{\sum }\dim (\pi )\left( \frac{\chi
_{\pi }(T)}{\dim (\pi )}\right) ^{\ell }\chi _{\pi }(g)}}\text{,}
\label{Slg-again}
\end{equation}%
and analyze each of the sub-sums $(S_{\ell ,g})_{k},$ $k=1,...,n,$ in (\ref%
{Slg-again}).

\begin{claim}
\label{C-Slgk}Suppose, $g\in (\partial (G))_{reg}$. Then, for $\ell \geq n,$

\begin{enumerate}
\item For $k=1$,$\ \ (S_{\ell ,g})_{1}=-\frac{2}{q}(\frac{1}{q})^{\ell
-n}+o(...);$

\item For $k=2$,\ $\ (S_{\ell ,g})_{2}=O(\frac{1}{q}(\frac{1}{q^{2}})^{\ell
-n});\vspace{0.01in}$

and more generally,

\item For $2\leq k\leq n$, $\ (S_{\ell ,g})_{k}=\left\{ 
\begin{array}{c}
O(\frac{1}{q^{\frac{k(k-1)}{2}}}(\frac{1}{q^{k}})^{\ell -n}),\text{ if }%
k\leq n-2; \\ 
O(\frac{1}{q^{\frac{n(n-1)}{2}}}(\frac{1}{q^{n-1}})^{\ell -n}),\text{ if }%
n-2<k.%
\end{array}%
\right. $
\end{enumerate}
\end{claim}

For a proof of Claim \ref{C-Slgk}, see Appendix \ref{P-C-Slgk}.\vspace{0.02in%
}

Finally, the estimates (\ref{Card-Mlg-again}) follows from Claim \ref{C-Slgk}%
. This completes the proof of Theorem \ref{T-STS}.

\section{\textbf{A Proof of the Agreement Conjecture\label{S-P-AC/EC}}}

In this section we propose a Fourier theoretic proof of the agreement
conjecture that uses the CRs estimates (\ref{CRs-GLn}), and certain curious
positivity results for the Fourier transform of the collection of matrices
of a fixed low enough rank.

\subsection{\textbf{The Statement}}

We know that for $k<\left\lfloor \frac{n}{2}\right\rfloor $, we have $(%
\widehat{GL}_{n})_{\otimes ,k}\subset (\widehat{GL}_{n})_{U,k}$.

We will show that the following is true:

\begin{theorem}[\textbf{Agreement}]
\label{T-Agreement}Suppose $k<\left\lfloor \frac{n}{2}\right\rfloor $. Then,
for sufficiently large $q$, we have,%
\begin{equation*}
(\widehat{GL}_{n})_{\otimes ,k}=(\widehat{GL}_{n})_{U,k}\text{.}
\end{equation*}
\end{theorem}

\subsection{\textbf{A Fourier Theoretic Input}}

We follow some of the development of Section \ref{S-UR}, in particular
Example \ref{E-FT}. Consider the collection $\mathcal{O}_{r}$ of rank $r$
matrices in $U\simeq M_{\left\lfloor \frac{n}{2}\right\rfloor
,(n-\left\lfloor \frac{n}{2}\right\rfloor )}$. \ Denote by $\widehat{1}_{%
\mathcal{O}_{r}}$ the Fourier transform of the characteristic function $1_{%
\mathcal{O}_{r}}$ of $\mathcal{O}_{r}$ (see Formula (\ref{Chi_rho_Or})). Let
us (abusing notation a bit) denote by $T\in U$, a conjugate of (\ref{T}),
i.e., a rank one matrix.

We have the following useful

\begin{fact}[Fourier transform of $\mathcal{O}_{r}$]
\label{F-FT-Or}For $r<\left\lfloor \frac{n}{2}\right\rfloor $, the value of $%
\widehat{1}_{\mathcal{O}_{r}}$ at $T$ is a positive integer, in fact,%
\begin{equation}
\widehat{1}_{\mathcal{O}_{r}}(T)=(\frac{1}{q^{r}}+o(...))\#(\mathcal{O}_{r}).
\label{FT-1Or-T}
\end{equation}
\end{fact}

Fact \ref{F-FT-Or} is a straightforward combination of Theorem \ref%
{T-FT-Mmnk-T}, and Part (1) of Lemma \ref{L-M,Mo-P1}, in Appendices \ref%
{A-Val-FT-Mmnk-T}, and \ref{P-P-FT-Mmnk-T}, respectively.

\subsection{\textbf{Proof of the Agreement Theorem}}

To prove Theorem \ref{T-Agreement}, take $\rho \in (\widehat{GL}_{n})_{U,k}$%
, $k<\left\lfloor \frac{n}{2}\right\rfloor ,$ and compute,%
\begin{eqnarray*}
\frac{\chi _{\rho }(T)}{\dim (\rho )} &=&\frac{\tsum\limits_{r=0}^{k}m_{r}%
\cdot \widehat{1}_{\mathcal{O}_{r}}(T)}{\tsum\limits_{r=0}^{k}m_{r}\cdot \#(%
\mathcal{O}_{r})} \\
&=&\frac{\tsum\limits_{r=0}^{k}m_{r}\cdot (\frac{1}{q^{r}}+o(...))\cdot \#(%
\mathcal{O}_{r})}{\tsum\limits_{r=0}^{k}m_{r}\cdot \#(\mathcal{O}_{r})}\text{
\ } \\
&\geq &(\frac{1}{q^{k}}+o(...))\frac{\tsum\limits_{r=0}^{k}m_{r}\cdot \#(%
\mathcal{O}_{r})}{\tsum\limits_{r=0}^{k}m_{r}\cdot \#(\mathcal{O}_{r})}=%
\frac{1}{q^{k}}+o(...),
\end{eqnarray*}%
where, the first equality is the expansion, of the restriction of $\rho $ to 
$U,$ discussed in Section \ref{S-UR} (see Formulas (\ref{rho_U}), (\ref%
{Chi_rho_Or}), and (\ref{rho-Or})), the second equality is by Formula (\ref%
{FT-1Or-T}), and finally, the inequality at the bottom is due to the
positivity of $\widehat{1}_{\mathcal{O}_{r}}(T)$ for $r<\left\lfloor \frac{n%
}{2}\right\rfloor $.

In particular, we see by the CRs estimates (\ref{CRs-GLn}) that $\rho $ must
be in $(\widehat{GL}_{n})_{\otimes ,k}$, as we wanted to show. This
completes the proof of Theorem \ref{T-Agreement}.

We proceed to the last section of this note, where we give some details on
how one might obtain the analytic results for the irreps of $GL_{n}$.

\section{\textbf{The eta Correspondence, the Philosophy of Cusp Forms, and
Tensor Rank\label{S-EC-TR}}}

We will wrap up the body of this note by giving some indications on how we
derived in \cite{Gurevich-Howe17, Gurevich-Howe19} the analytic results
described in Section \ref{S-AI-TRk-GLn}. In this way or another, this means
to address the following:\smallskip

\textbf{Question: }How to get information on the (e.g., CRs, dimensions, and
cardinality of the set of) $\otimes $-rank $k$ irreps of $GL_{n}$?\smallskip

One way to answer this question was carried out in \cite{Gurevich-Howe17,
Gurevich-Howe19}. It used the philosophy of cusp forms, and developed
criteria for representations to be of tensor rank $k$ \ in terms of their
appearance in representations induced from parabolic subgroups.

In more detail, the process of getting the information on the tensor rank $k$
irreps includes the following three steps:\vspace{0.01in}

\begin{enumerate}
\item \textbf{Eta correspondence (EC).} To some extent the EC might be
considered as giving you a convenient place where to search for a formula
for the irreps of tensor rank $k$. Moreover, it allows one to count the
number of such irreps.\vspace{0.01in}

\item \textbf{Philosophy of cusp forms (P-of-CF). }This is a method, put
forward in the 60s by Harish-Chandra \cite{Harish-Chandra70}, that allows
one to write formulas for irreps of groups like $GL_{n}$. In particular, in
\ our case the EC from Step (1) above lead us to find certain P-of-CF
formulas that seems to be effective for the analysis we want to do for the
irreps of tensor rank $k$.

\item \textbf{Derivation of the analytic information. }Having the formulas
from Step (2) above, one, in principal, does explicit calculations and
derives the analytic results.
\end{enumerate}

Let us go over the main statement of Step (1) above, then write down the
P-of-CF formulas of Step (2) as they applied to tensor rank $k$ irreps, and
finally give the main computations done in Step (3) in order to derive: the
CRs at the transvection and dimensions for the $\otimes $-rank $k$ irreps of 
$GL_{n}$, and the cardinality of the set of all these irreps.

\subsection{\textbf{The eta Correspondence and Strict Tensor Rank\label%
{S-EC-STR}}}

Recall (see Section \ref{S-TR}, in particular Definition \ref{D-TR}) that an
irrep of $GL_{n}$ is of tensor rank $k,$ $0\leq k\leq n$, if up to twist by
a character (one dim irrep) it appears in $\omega _{n,k}=\omega
_{n}^{\otimes ^{k}}$and not in $\omega _{n,(k-1)}$, where $\omega _{n,k}$
denotes the permutation representation of $GL_{n}$ on the space $%
L^{2}(M_{n,k})$. Let us introduce the following terminology:

\begin{definition}
\label{D-STR}We say that an irrep $\rho $ of $GL_{n}$ is of \underline{%
strict tensor rank} $k$, $0\leq k\leq n$, if it appears in $\omega _{n,k}$,
but not in $\omega _{n,(k-1)}$.
\end{definition}

Let us denote the set of all irreps of $GL_{n}$ of strict tensor rank $k$ by 
$(\widehat{GL}_{n})_{\otimes ,k}^{\star }$.

Since every irrep of tensor rank $k$ is up to twist by a character in $(%
\widehat{GL}_{n})_{\otimes ,k}^{\star }$, and this twist does not affect the
dimension or the CR of the transvection, so it might be beneficial for us to
get information on the members of $(\widehat{GL}_{n})_{\otimes ,k}^{\star }$%
. Since they all appear inside $L^{2}(M_{n,k}),$ we want to zoom into this
space and locate them. To do this, in \cite{Gurevich-Howe17, Gurevich-Howe19}
we followed \cite{Howe73} and (as in Section \ref{S-T-EC}, but now for any
value of $k$) use the concept of $GL_{n}$-$GL_{k}$ dual pair.

Consider the oscillator rep $\omega _{n,k}$ as the joint action (\ref%
{omega_n,k}) of $GL_{n}\times GL_{k}$ on $L^{2}(M_{n,k})$, and decompose it
as in (\ref{Isot-Dec}) to a direct sum of $GL_{k}$-isotypic components%
\begin{equation}
\omega _{n,k}\simeq \sum_{\tau \in \widehat{GL}_{k}}\mathcal{M}(\tau
)\otimes \tau ,  \label{omega-nk-again}
\end{equation}%
where each multiplicity space $\mathcal{M}(\tau )$ is a rep of $GL_{n}$.

It turns out that (similar to the $k<\left\lfloor \frac{n}{2}\right\rfloor $
case discussed in Section \ref{S-T-EC}) for "most" $\tau $'s the space $%
\mathcal{M}(\tau )$ contains a distinguished irrep $\eta (\tau )$, which is
in fact from $(\widehat{GL}_{n})_{\otimes ,k}^{\star }$. To describe it more
closely, let us consider the parabolic subgroup $P_{k,n-k}\subset GL_{n}$
stabilizing the first $k$ coordinates subspace of $\mathbb{F}_{q}^{n}$ (we
denoted it by $P_{k}$ previously, see Formulas (\ref{Xm}) and (\ref{Pm})),
and recall that it has a natural projection onto its Levi component $%
GL_{k}\times GL_{n-k}$. Then, to each irrep $\tau \in \widehat{GL}_{k}$ we
can consider the rep $\tau \otimes \mathbf{1}_{n-k}$ of $GL_{k}\times
GL_{n-k},$ pull it back to $P_{k,n-k}$ and look at the induced representation%
\begin{equation}
I_{\tau }=Ind_{P_{k,n-k}}^{GL_{n}}(\tau \otimes \mathbf{1}_{n-k}).
\label{I-tau}
\end{equation}

Now we can write down, with some more details, the natural extension given
in \cite{Gurevich-Howe17, Gurevich-Howe19} for the eta correspondence
described in Theorem \ref{T-EC-U}.

\begin{theorem}[\textbf{eta correspondence}]
\smallskip \label{T-EC}Take $\tau \in \widehat{GL}_{k},$ $0\leq k\leq n,$
and look at the decomposition (\ref{omega-nk-again}). We have,

\begin{enumerate}
\item \textbf{Existence. }The representation $\mathcal{M}(\tau )$ contains a
strict tensor rank $k$ component if and only if $\tau $ is of strict tensor
rank \ $\geq 2k-n.$

Moreover, if the condition of Part (1) is satisfied, then,

\item \textbf{Uniqueness.}\textit{\ } the representation $\mathcal{M}(\tau )$
has a unique constituent $\eta (\tau )$ of strict tensor rank $k$, and it
appears with multiplicity one.\smallskip

and,\smallskip

\item \textbf{Approximate formula. }the constituent $\eta (\tau )$ satisfies 
$\eta (\tau )<I_{\tau }<\mathcal{M}(\tau )$, and we have,%
\begin{equation*}
\ I_{\tau }=\eta (\tau )+\dsum\limits_{\rho }\rho ,\ 
\end{equation*}%
where the sum is multiplicity free, and over certain irreps $\rho $ which
are of strict tensor rank less than $k$ and dimension smaller than $\eta
(\tau )$.

Finally, the mapping 
\begin{equation}
\tau \longmapsto \eta (\tau )\text{,}  \label{etaC-STR}
\end{equation}%
gives an explicit bijective correspondence 
\begin{equation*}
(\widehat{GL}_{k})_{\otimes ,\geq 2k-n}^{\star }\longrightarrow (\widehat{GL}%
_{n})_{\otimes ,k}^{\star },
\end{equation*}%
between the collection $(\widehat{GL}_{k})_{\otimes ,\geq 2k-n}^{\star }$ of
irreps of $GL_{k}$ of strict tensor rank $\geq 2k-n$, and the set $(\widehat{%
GL}_{n})_{\otimes ,k}^{\star }$ of strict tensor rank\textit{\ }$k$ irreps
of $GL_{n}$.
\end{enumerate}
\end{theorem}

Next, we want to analyze further the $\eta (\tau )$'s mentioned just above,
and find more about them.

\subsection{\textbf{The Philosophy of Cusp Forms Formula and Rank\label%
{S-PofCF-F}}}

Using the eta correspondence and in particular the observation that each of
the $\eta (\tau )$ (\ref{etaC-STR}) is by Part (3) of Theorem \ref{T-EC}
"not too far" from being the induced representation $I_{\tau }$ (\ref{I-tau}%
), we were lead in \cite{Gurevich-Howe17, Gurevich-Howe19} to an explicit
Harish-Chandra type formula for irreps of $GL_{n}$ from which one can easily
read off their strict tensor rank and tensor rank.

In this note we will just write down the above mentioned formula and explain
how to get the rank invariants from it. We will leave the details of how we
arrived to that expression (as well as the details of the relevant
Harish-Chandra P-of-CF theory) to \cite{Gurevich-Howe17, Gurevich-Howe19}.

\subsubsection{\textbf{The P-of-CF Formula}}

The formula mentioned just above will be given in terms of representations
induced from certain representations of parabolic subgroups of $GL_{n}$ that
contain the standard Borel subgroup $B$ of upper-triangular matrices \cite%
{Borel69}.

Recall that to every ordered partition 
\begin{equation}
D=\{d_{1}\geq ...\geq d_{\ell }\},  \label{D}
\end{equation}%
of $n$, we can associate the \textit{standard flag }$F_{D}$ of subspaces of $%
\mathbb{F}_{q}^{n}$:%
\begin{equation}
F_{D}:\text{ \ }0=X_{m_{0}}\subset X_{m_{1}}\subset ...\subset X_{m_{\ell }}=%
\mathbb{F}_{q}^{n},  \label{FD}
\end{equation}%
where for each $1\leq j\leq \ell $, we have $m_{j}-m_{j-1}=d_{j}$, and $%
X_{m_{j}}$ is the first $m_{j}$-coordinates subspace of $\mathbb{F}_{q}^{n}$
(see Formula (\ref{Xm})).

In particular, we can attach to $D$ (\ref{D}) the parabolic parabolic
subgroup%
\begin{equation}
P_{D}=Stab_{GL_{n}}(F_{D}),  \label{PD}
\end{equation}%
of all elements $g\in GL_{n}$ that stabilize the flag $F_{D}$ (\ref{FD}),
i.e., satisfy $g(X_{m_{j}})=X_{m_{j}}$, for all $j$.

Note that the subgroup $P_{D}$ has the following structure of a block upper
triangular matrix:%
\begin{equation*}
P_{D}=\left\{ \left( 
\begin{tabular}{cccc}
\cline{1-1}
\multicolumn{1}{|c}{$C_{1}$} & \multicolumn{1}{|c}{$\ast $} & $\ast $ & $%
\ast $ \\ \cline{1-2}
& \multicolumn{1}{|c}{$\ddots $} & \multicolumn{1}{|c}{$\ast $} & $\ast $ \\ 
\cline{2-2}
&  & $\ddots $ & $\ast $ \\ \cline{4-4}
&  &  & \multicolumn{1}{|c|}{$C_{\ell }$} \\ \cline{4-4}
\end{tabular}%
\right) ;\text{ }C_{j}\in GL_{d_{j}},\text{ }j=1,...,\ell \right\} ,
\end{equation*}%
and in particular admits a natural projection%
\begin{equation}
\overset{P_{D}}{\overbrace{\left\{ \left( 
\begin{tabular}{cccc}
\cline{1-1}
\multicolumn{1}{|c}{$C_{1}$} & \multicolumn{1}{|c}{$\ast $} & $\ast $ & $%
\ast $ \\ \cline{1-2}
& \multicolumn{1}{|c}{$\ddots $} & \multicolumn{1}{|c}{$\ast $} & $\ast $ \\ 
\cline{2-2}
&  & $\ddots $ & $\ast $ \\ \cline{4-4}
&  &  & \multicolumn{1}{|c|}{$C_{\ell }$} \\ \cline{4-4}
\end{tabular}%
\right) \right\} }}\longrightarrow \overset{L_{D}}{\overbrace{\left\{ \left( 
\begin{tabular}{cccc}
\cline{1-1}
\multicolumn{1}{|c}{$C_{1}$} & \multicolumn{1}{|c}{} &  &  \\ \cline{1-2}
& \multicolumn{1}{|c}{$\ddots $} & \multicolumn{1}{|c}{} &  \\ \cline{2-2}
&  & $\ddots $ &  \\ \cline{4-4}
&  &  & \multicolumn{1}{|c|}{$C_{\ell }$} \\ \cline{4-4}
\end{tabular}%
\right) \right\} }},  \label{LD}
\end{equation}%
onto its \textit{Levi component }$L_{D}\simeq GL_{d_{1}}\times ...\times
GL_{d_{\ell }}$.\vspace{0.01in}

\paragraph{\textbf{Split Representations}}

The first type of representations we will need for our formula have been
called in \cite{Gurevich-Howe17, Gurevich-Howe19} \textit{split principal
series. }They are the constituents of the induced representations $%
Ind_{B}^{GL_{n}}(\mathbf{\chi )}$ from characters of the Borel subgroup $B$.
To write efficient expressions for them, let us first look at the \textit{%
spherical principal series (SPS) }representations, which are the
constituents of the induced representation $Ind_{B}^{GL_{n}}(\mathbf{1)}$
from the trivial representation $\mathbf{1}$ of $B$.

The SPS irreps can be realized nicely using the parabolic subgroups $P_{D}$ (%
\ref{PD}). Indeed, consider the induced representation%
\begin{equation}
I_{D}=Ind_{P_{D}}^{GL_{n}}(\mathbf{1).}  \label{ID}
\end{equation}%
Recall that,

\begin{definition}
If, in addition to $D$ (\ref{D}), we have another partition $D^{\prime
}=\{d_{1}^{\prime }\geq d_{2}^{\prime }\geq ...\geq d_{\ell ^{\prime
}}^{\prime }\}$ of $n,$ then we say that $D^{\prime }$ \underline{dominates} 
$D$, and write $D\preceq D^{\prime },$ if $\ell ^{\prime }\leq \ell $ and%
\begin{equation*}
\sum_{j=1}^{j}d_{j}\leq \sum_{j=1}^{j}d_{j}^{\prime },\text{ \ \ \ for }%
j=1,...,\ell ^{\prime }.
\end{equation*}
\end{definition}

With the above terminology one can show \cite{Gurevich-Howe17,
Gurevich-Howe19},

\begin{fact}
The representation $I_{D}$ contains a constituent 
\begin{equation}
\rho _{D}<I_{D},  \label{rhoD}
\end{equation}
with multiplicity one, and with the property that it is not contained in any 
$I_{D^{\prime }}$ with $D^{\prime }\succneqq D$ in the dominance order.
\end{fact}

\begin{remark}
The representation $\rho _{D}$ can also be distinguished by its dimension:
it is the only constituent of $I_{D}$ whose dimension, as a polynomial in $q$%
, has the same degree as the cardinality of $GL_{n}/P_{D}$.
\end{remark}

The irreps $\rho _{D}$ (\ref{rhoD}), where $D$ runs over all ordered
partitions of $n$, are pairwise non-isomorphic and exhaust the collection of
SPS irreps.

Using the $\rho _{D}$'s we can describe the split principal series
representations. Indeed, consider the parabolic subgroup $P_{S}$ attached to
the partition $S=\{s_{1}\geq ...\geq s_{\ell }\}$ of $n$, and its Levi
component $L_{S}\simeq GL_{s_{1}}\times ...\times GL_{s_{\ell }}$. Suppose $%
D_{1},...,D_{\ell },$ are ordered partitions of $s_{1},...,s_{\ell },$
respectively, and that $\chi _{1},...,\chi _{\ell }$, are $\ell $ distinct
characters of $\mathbb{F}_{q}^{\ast }$. Then we can consider the
representation $\tbigotimes\limits_{j=1}^{\ell }[(\chi _{j}\circ \det
)\otimes \rho _{D_{j}}]$ of $L_{S}$, pull it back to $P_{S}$, and induce to
form,%
\begin{equation}
\rho _{S}=Ind_{P_{S}}^{GL_{n}}(\tbigotimes\limits_{j=1}^{\ell }[(\chi
_{j}\circ \det )\otimes \rho _{D_{j}}]).  \label{rhoS}
\end{equation}%
It is not difficult to check that (up to order of the inducing factors that
correspond to $s_{j}$'s of the same size) the irreps $\rho _{S}$'s \ref{rhoS}
are irreducible, pairwise non-isomorphic, and exhausts the collection of
split principal series representations.\vspace{0.01in}

\paragraph{\textbf{Unsplit Representations}}

The second type of representations we use in our formula are the irreps we
called in \cite{Gurevich-Howe17} \textit{unsplit. }To define them, let us
first recall the following basic objects in Harish-Chandra's P-of-CF \cite%
{Harish-Chandra70}:

\begin{definition}
A representation $\kappa $ of $GL_{n}$ is called \underline{\textit{cuspidal}%
} if it does not contain a non-trivial fixed vector for the unipotent
radical of any parabolic subgroup stabilizing a flag in $\mathbb{F}_{q}^{n}$.
\end{definition}

In the above definition, it is enough to consider parabolic subgroups of the
form $P_{D}$ (\ref{PD}), and their unipotent radicals, i.e., the kernels of
the projections (\ref{LD}).

In this note we will not explicitly discuss the cuspidal representations
(for this see \cite{Gel'fand70, Howe-Moy86, Zelevinsky81}), but only use
them and some of their properties as needed. In particular, if $%
U=\{u_{1}\geq ...\geq u_{m}\}$ is an ordered partition of $n,$ and $P_{U}$
is the corresponding parabolic subgroup, with Levi component $L_{U}\simeq
GL_{u_{1}}\times ...\times GL_{u_{m}}$, we can take cuspidal irreps $\kappa
_{1},...,\kappa _{m}$, of $GL_{u_{1}},...,GL_{u_{m}},$ respectively, then
form the representation $\tbigotimes\limits_{i=1}^{m}\kappa _{i}$ of $L_{U%
\text{,}}$ and induce to get $Ind_{P_{U}}^{GL_{n}}(\tbigotimes%
\limits_{i=1}^{m}\kappa _{i}).$

If in the partition $U$ above, we have that $u_{i}\geq 2,$ for every $%
i=1,...,m,$ and 
\begin{equation}
\rho _{U}<Ind_{P_{U}}^{GL_{n}}(\tbigotimes\limits_{i=1}^{m}\kappa _{i}),
\label{rhoU}
\end{equation}%
is an irreducible component, then we call $\rho _{U}$ \textit{unsplit
representation.}

\paragraph{\textbf{General Representations}}

The P-of-CF formula for general irreps of $GL_{n},$ is obtained by a
combination of the formulas for split and unsplit representations defined
just above.

Suppose $n=u+s,$ with integers $u,s\geq 0$. Then we can consider the
parabolic subgroup $P_{u,s}\subset GL_{n}$, stabilizing the first $u$%
-coordinates subspace of $\mathbb{F}_{q}^{n}$, and its Levi subgroup $%
L_{u,s}\simeq GL_{u}\times GL_{s}$. Take an unsplit irrep $\rho _{U}$ of $%
GL_{u}$ in addition to a split irrep $\rho _{S}$ of $GL_{s},$ then pullback
the rep $\rho _{U}\otimes \rho _{S}$ of $L_{u,s}$ to $P_{u,s}$ and induce to
get%
\begin{equation}
\rho _{U,S}=Ind_{P_{u,s}}^{GL_{n}}(\rho _{U}\otimes \rho _{S}).
\label{rho_US-1}
\end{equation}%
Now, the philosophy of cusp forms \cite{Bump04, Howe-Moy86, Harish-Chandra70}
tells us that\vspace{0.01in}

\textbf{(a)} $\rho _{U,S}$ is irreducible; and\vspace{0.01in}

\textbf{(b)} the map $(\rho _{U},\rho _{S})\mapsto \rho _{U,S}$ is an
injection from the relevant subsets of the unitary dual of $GL_{u}\times
GL_{s}$ into the unitary dual of $GL_{n}$; and\vspace{0.01in}

\textbf{(c)} all irreducible representations of $GL_{n}$ arise in this way
(including the cuspidal representations, which are included in the situation
when $P_{n,0}=GL_{n}$).

For a later use, we want to make Formula (\ref{rho_US-1}) a bit more
explicit. Indeed, suppose, in addition, that our $S=\{s_{1}\geq ...\geq
s_{\ell }\}$ is an ordered partition of $s,$ then we can consider the
standard parabolic subgroup $P_{u,s_{1},...,s_{\ell }}\subset GL_{n}$ with
blocks of sizes $u,s_{1},...,s_{\ell },$ and pullback to it the rep (with
distinct $\chi _{j}$'s) $\rho _{U}\otimes \lbrack
\tbigotimes\limits_{j=1}^{\ell }((\chi _{j}\circ \det )\otimes \rho
_{D_{j}})]$ of its Levi component $L_{u,s_{1},...,s_{\ell }}\simeq
GL_{u}\times GL_{s_{1}}\times ...\times GL_{s_{\ell }},$ and induce to get 
\begin{equation}
\rho _{U,S}=Ind_{P_{u,s_{1},...,s_{\ell }}}^{GL_{n}}\left( \rho _{U}\otimes 
\left[ \tbigotimes\limits_{j=1}^{\ell }(\chi _{j}\circ \det )\otimes \rho
_{D_{j}}\right] \right) \text{.}  \label{rho_US-2}
\end{equation}

We will call (\ref{rho_US-2}) the \textbf{P-of-CF formula.}

\subsubsection{\textbf{Reading Ranks from the P-of-CF Formula}}

In \cite{Gurevich-Howe17, Gurevich-Howe19} we showed that one can compute
the strict tensor rank and tensor rank of a representation from its P-of-CF
Formula (\ref{rho_US-2}), more precisely directly from its split principal
series component. To state this, and similar results, it is convenient to
use the notions of \underline{\textit{tensor co-rank}} and \underline{%
\textit{strict tensor co-rank}}, by which we mean, respectively, $n$ minus
the tensor rank and $n$ minus the strict tensor rank.

\begin{fact}
\label{F-TR-from-PofCF-F}We have,

\begin{enumerate}
\item For an ordered partition $D=\{d_{1}\geq ...\}$ of $n$, the tensor
co-rank of the SPS representation $\rho _{D}$ (\ref{PD}) is the same as its
strict tensor co-rank and is equal to $d_{1}$.

\item The tensor co-rank of the representation $\rho _{U,S}$ of $GL_{n}$
described by Formula (\ref{rho_US-2}), is the maximum of the tensor co-ranks
of the SPS representations $\rho _{D_{j}},$ $j=1,...,\ell $, that appear in
description of the split part of $\rho _{U,S}$. The strict tensor co-rank of 
$\rho _{U,S}$ is the strict tensor rank of the SPS representation $\rho
_{D_{j}}$ that is twisted in (\ref{rho_US-2}) by the trivial character.
\end{enumerate}
\end{fact}

\subsection{\textbf{Deriving the Analytic Information for Tensor Rank }$k$%
\textbf{\ irreps of }$GL_{n}$}

In this last section we want to remark briefly on how the eta correspondence
(Section \ref{S-EC-STR}) and the P-of-CF Formula for tensor rank $k$ irreps
of $GL_{n}$, enable us in \cite{Gurevich-Howe19} to obtain the analytic
information for these irreps given in Section \ref{S-AI-TRk-GLn}.

\subsubsection{\textbf{Character Ratios on The Transvection}}

It is not difficult to see \cite{Gurevich-Howe19} that one just need to
estimate the CR at $T$ (\ref{T}) for irreps of tensor rank $k$ of the form

\begin{equation}
\rho _{U,D}=Ind_{P_{u,d}}^{GL_{n}}\left( \rho _{U}\otimes \rho _{D}\right) ,
\label{rho_UD}
\end{equation}%
where

\begin{itemize}
\item $P_{u,d}$ is the parabolic subgroup with blocks of sizes $u$ and $d$, $%
u+d=n$;

\item $\rho _{U}$ is an unsplit irrep (see (\ref{rhoU})) of $GL_{u}$;

\item $D=\{d_{1}\geq ...\}$ is a partition of $d$, with longest row of
length $d_{1}=n-k;\vspace{0.01in}$

and

\item $\rho _{D}$ the SPS irreps (\ref{PD}) attached to $D$.
\end{itemize}

In particular, using the standard formula \cite{Fulton-Harris91} for
character of induced representation, the estimate is reduced in \cite%
{Gurevich-Howe19} to the known cardinality of $GL_{n}/P_{u,d}$, and the
following two specific cases of CR estimates:\vspace{0.02in}

\textbf{SPS case. }Consider a SPS irrep $\rho _{D}$ (\ref{PD}) of $GL_{n},$
where $D=\{d_{1}\geq ...\geq d_{\ell }\}$ is a partition of $n$. Then, the
tensor rank of $\rho _{D}$ is equal to $k=n-d_{1}$ (see Fact \ref%
{F-TR-from-PofCF-F}), and,%
\begin{equation*}
\frac{\chi _{\rho _{D}}(T)}{\dim (\rho _{D})}=\left\{ 
\begin{array}{c}
\frac{1}{q^{n-d_{1}}}+o(...)\text{, \ \ if \ }d_{1}>d_{2}\text{;} \\ 
\\ 
\frac{c_{D}}{q^{n-d_{1}}}+o(...)\text{, \ \ \ \ otherwise,}%
\end{array}%
\right. \text{,}
\end{equation*}%
where $c_{D}$ is a certain integer depending only on $D$ (and not on $q$).%
\vspace{0.02in}

\textbf{Unsplit case. }Consider an unsplit irrep $\rho _{U}$ (\ref{rhoU}) of 
$GL_{n}$. Then, the tensor rank of $\rho _{U}$ is equal to $n$, and 
\begin{equation*}
\frac{\chi _{\rho _{U}}(T)}{\dim (\rho _{U})}=\frac{-1}{q^{n-1}-1}.
\end{equation*}

\subsubsection{\textbf{Dimensions of Irreps\label{S-Dim-TR-k}}}

Again, it is enough to compute the dimensions of irreps of the form $\rho
_{U,D}$ (\ref{rho_UD}), and for such irrep it boils down to the following
two known \cite{Gurevich-Howe19} cases:\vspace{0.02in}

\textbf{SPS case. }The dimension of the SPS irrep $\rho _{D}$ (\ref{rhoD})
attached to a partition $D=\{d_{1}\geq d_{2}\geq ...\geq d_{\ell }\}$ of $n,$
satisfies,%
\begin{equation*}
\dim (\rho _{D})=q^{d_{D}}+o(...),
\end{equation*}%
where $d_{D}=\tsum\limits_{1\leq i<j\leq \ell }d_{i}d_{j}$.\vspace{0.02in}

\textbf{Cuspidal case. }The dimension of any cuspidal representation of $%
GL_{n}$ is $q^{\frac{n(n-1)}{2}}+o(...)$ \cite{Bump04, Gel'fand70}.

\subsubsection{\textbf{The Number of Irreps of Tensor Rank }$k$}

The P-of-CF formula for the SPS irreps shows that the cardinality of $(%
\widehat{GL}_{n})_{\otimes ,n-1}$ is $c_{n-1}q^{n}+o(...)$, for some
constant $c_{n-1}>0$, independent of $q$, and the well known \cite{Bump04,
Gel'fand70} cardinality of the collection of cuspidal irreps, implies that $%
\#((\widehat{GL}_{n})_{\otimes ,n})$ is also of the form $c_{n}q^{n}+o(...)$%
, for some constant $c_{n}>0$, independent of $q$, and using more or less
the same knowledge one can show that $c_{n-1}+c_{n}\geq 1,$ and hence $=1$,
since the size of $\widehat{GL}_{n}$ is $q^{n}+o(...).$

For $k<n-1$, the estimation given in \cite{Gurevich-Howe19} use the eta
correspondence\textbf{\ }described by Theorem \ref{T-EC}. For such $k,$ the
domain of the eta map (\ref{etaC-STR})%
\begin{equation*}
\eta :(\widehat{GL}_{k})_{\otimes ,\geq 2k-n}^{\star }\longrightarrow (%
\widehat{GL}_{n})_{\otimes ,k}^{\star },
\end{equation*}%
has cardinality $q^{k}+o(...),$ for example because it contains the the
tensor rank $k-1,$ and tensor rank $k$ irreps of $GL_{k}$. Moreover, using
the P-of-CF formula one can show (see \cite{Gurevich-Howe19}) that for $%
q^{k}+o(...)$ of the irreps $\tau \in (\widehat{GL}_{k})_{\otimes ,\geq
2k-n}^{\star }$, the irrep $\eta (\tau )$ is of tensor rank $k,$ and stays
like this even after twists by the $q-1$ characters $GL_{n}$, and, moreover,
all these are pairwise non-isomorphic. This shows that $\#((\widehat{GL}%
_{n})_{\otimes ,k})=q^{k+1}+o(...),$ for $k<n-1$, as claimed.

This completes the story we were trying to give in the body of this note.

\appendix

\section{\textbf{Fourier Transform of Sets of Matrices of Fixed Rank}\label%
{A-FT}}

For $k\leq m\leq n$, we look at the Fourier transform (FT) of the set $%
(M_{m,n})_{k},$ of $m\times n$ matrices of rank $k$ over $\mathbb{F}_{q}$,
and want to evaluate it on a rank one matrix.\smallskip

Let us identify, in the standard way, the space $M_{m,n}$ of $m\times n$
matrices over $\mathbb{F}_{q}$ with its vector space dual via the trace map,
and then to its Pontryagin dual in the standard way, after fixing an
additive character $\psi $ of $\mathbb{F}_{q}.$ In particular, we can
associate with a function $f$ on $M_{m,n}$ its Fourier transform $\widehat{f}
$ given by 
\begin{equation}
\widehat{f}(B)=\sum_{A\in M_{m,n}}f(A)\psi (-trace(B^{t}\circ A)),\text{ \ }%
B\in M_{m,n}\text{.}  \label{FT}
\end{equation}

\subsection{\textbf{A\ Formula for the Fourier Transform of }$(M_{m,n})_{k}$}

By the FT of $(M_{m,n})_{k},$ we mean the FT of the characteristic function $%
1_{(M_{m,n})_{k}}$ of that set. Using (\ref{FT}), it is given by%
\begin{equation}
\widehat{1}_{(M_{m,n})_{k}}(B)=\sum_{A\in (M_{m,n})_{k}}\psi
(trace(B^{t}\circ A)),\text{ \ }B\in M_{m,n}\text{.}  \label{FT-Mmnk}
\end{equation}

In particular, $\widehat{1}_{(M_{m,n})_{k}}(B)$ is a real number.

Next, denote by $(M_{m,n})_{k}^{o}(B)$ the set of matrices $A\in
(M_{m,n})_{k}$ such that $trace(B^{t}\circ A)=0$.

\begin{claim}
\label{C-F-FT-Mmnk}We have,%
\begin{equation}
\widehat{1}_{(M_{m,n})_{k}}(B)=-(\frac{1}{q-1})\#((M_{m,n})_{k})+(\frac{q}{%
q-1})\#((M_{m,n})_{k}^{o}(B))\text{.}  \label{F-FT-Mmnk}
\end{equation}
\end{claim}

For the calculations leading to Formula (\ref{F-FT-Mmnk}), see Appendix \ref%
{P-C-F-FT-Mmnk}.

Note that, for a non-zero $B\in M_{m,n}$, the map $A\mapsto trace(B^{t}\circ
A)$ is a linear functional on the space $M_{m,n}$, so its kernel is a
hyperplane, i.e., it has co-dimension $1$. So $(M_{m,n})_{k}^{o}(B)$ is some
sort of hypersurface in $(M_{m,n})_{k}$, and it should be reasonable to
compare the two numbers appearing as the two terms of the sum on the right
side of (\ref{F-FT-Mmnk}), and see if it is negative or positive.

We proceed to do just this in the case of a rank one matrix.

\subsection{\textbf{The Value of the FT of }$(M_{m,n})_{k}$\textbf{\ on a
Rank One Matrix\label{A-Val-FT-Mmnk-T}}}

Let $T\in M_{m,n}$ be a rank one matrix. The value of $\widehat{1}%
_{(M_{m,n})_{k}}$ at $T$ is of course independent of $T$, and it can be
computed explicitly. To write it down and use it, it will be useful for us
to recall (see also \cite{Artin57}) that the cardinality of the group $%
GL_{k}=GL_{k}(\mathbb{F}_{q})$ is 
\begin{eqnarray*}
\#(GL_{k}) &=&(q^{k}-1)(q^{k}-q)(q^{n}-q^{2})\cdots (q^{n}-q^{n-1}) \\
&=&q^{\frac{k(k-1)}{2}}(\prod\limits_{a=1}^{k}(q^{a}-1)) \\
&=&q^{n^{2}}+o(...),
\end{eqnarray*}%
and that the cardinality of $\Gamma _{n,k}$, the Grassmannian of all $k$%
-dimensional subspace of $\mathbb{F}_{q}^{n}$, is%
\begin{eqnarray}
\#(\Gamma _{n,k}) &=&\frac{\#(GL_{n})}{q^{k(n-k)}\cdot \#(GL_{k})\cdot
\#(GL_{n-k})}  \label{Card-Gama-nk} \\
&=&\frac{\tprod\nolimits_{a=1}^{k}(q^{n-k+a}-1)}{\tprod%
\nolimits_{a=1}^{k}(q^{a}-1)}  \notag \\
&=&q^{k(n-k)}+o(...).  \notag
\end{eqnarray}

\begin{remark}
In order for certain formulas to include all cases, we may sometime use the
notation $\Gamma _{n,k}$ also in the case $k>n$, by this we mean the empty
set with $\#(\Gamma _{n,k})=0.$
\end{remark}

Now we can write an explicit expression,

\begin{theorem}[\textbf{FT of }$M_{(m,n);k}$\textbf{\ on rank one matrix}]
\label{T-FT-Mmnk-T}Assume $k\leq m\leq n$. Then, the value of the FT of $%
(M_{m,n})_{k}$ at a rank one matrix $T\in M_{m,n},$ is 
\begin{equation*}
\widehat{1}_{(M_{m,n})_{k}}(T)=(q^{n+m-k}-q^{n}-q^{m}+1)(\frac{\#(\Gamma
_{n,k})\#(\Gamma _{m,k})\#(GL_{k})}{(q^{n}-1)(q^{m}-1)}).
\end{equation*}%
In particular, it is positive if $k<m$, and negative if $k=m.$
\end{theorem}

For a proof of Theorem \ref{T-FT-Mmnk-T}, see Appendix \ref{P-P-FT-Mmnk-T}.

\section{\textbf{Proofs}}

\subsection{\textbf{Proofs for Section \protect\ref{S-In}\label{Pr-S-In}}}

\subsubsection{\textbf{Proof of Proposition \protect\ref{P-C-Mlg}\label%
{Pr-P-C-Mlg}}}

\begin{proof}
For the basic notions and facts from representation theory of finite groups
see \cite{Serre77}.

The result we need to prove holds for any finite group $G$ and any conjugacy
class $C\subset G$.

For a finite group $G$, consider the regular representation $\pi _{G}$ on
the space $L(G)$ of complex-valued functions on $G$.

Associated with the conjugacy class $C\subset G$, we have a summation
operator $A_{C}:L(G)\rightarrow L(G),$ given by 
\begin{equation*}
A_{C}=\sum_{c\in C}\pi _{G}(c)\text{.}
\end{equation*}%
It is easy to check that%
\begin{equation*}
\#(M_{\ell ,g})=\left[ \left( A_{C}\right) ^{\ell }\delta _{1}\right] (g)%
\text{,}
\end{equation*}%
where $\delta _{1}$ denotes the Dirac delta function at the identity element
of $G$.

Next, for each irrep $\pi \in \widehat{G}$, denote by $Pr_{\pi }$ the
projection of $L(G)$ onto the $\pi $-isotypic component. Then we
have:\smallskip

(*) $\ \sum\limits_{\pi \in \widehat{G}}Pr_{\pi }=Id,$ the identity operator
on $L(G)$;\smallskip

and,\smallskip

(**) $\ \left( A_{C}\right) ^{\ell }\circ Pr_{\pi }=\#(C^{\ell })\cdot
\left( \frac{\chi _{\pi }(C)}{\dim (\pi )}\right) ^{\ell }\cdot Pr_{\pi }$%
.\smallskip

Fact (**) above follows by a direct computation using Schur's lemma, that
holds here since $A_{C}$ intertwines the action of $G$ on each $\pi \in 
\widehat{G}$.

Now, the projectors $Pr_{\pi }$, $\pi \in \widehat{G}$, have an explicit
formula \cite{Serre77},%
\begin{equation}
Pr_{\pi }=\frac{\dim (\pi )}{\#(G)}\cdot \sum_{h\in G}\chi _{\pi
}(h^{-1})\pi _{G}(h)\text{.}  \label{Pr_pi}
\end{equation}%
Concluding, we get,%
\begin{eqnarray*}
\#(M_{\ell ,g}) &=&\left[ \left( A_{C}\right) ^{\ell }\delta _{1}\right] (g)
\\
&=&\sum\limits_{\pi \in \widehat{G}}\left[ \left( A_{C}\right) ^{\ell }\circ
Pr_{\pi }(\delta _{1})\right] (g) \\
&=&\frac{\#(C^{\ell })}{\#(G)}\left( \sum\limits_{\pi \in \widehat{G}}\dim
(\pi )\left( \frac{\chi _{\pi }(C)}{\dim (\pi )}\right) ^{\ell }\chi _{\pi
}(g^{-1})\right) ,
\end{eqnarray*}%
where for the second equality we used Fact (*) above, and the third equality
is obtained using Fact (**), and Identity (\ref{Pr_pi}). This completes the
derivation of Formula (\ref{F-C-Mlg}), as needed.
\end{proof}

\subsection{\textbf{Proofs for Section \protect\ref{S-CRs-URank}}}

\subsubsection{\textbf{Proof of Proposition} \protect\ref{P-mult-rho_U} 
\label{P-P-mult-rho_U}}

\begin{proof}
Take a rep $\rho $ of $GL_{n}$. Then for every element $l$ in the Levi
subgroup $L_{\left\lfloor \frac{n}{2}\right\rfloor }$ (\ref{Levi-decomp}),
we have the representation $\rho ^{l}$ of $GL_{n},$ given by $\rho
^{l}(g)=\rho (l\cdot g\cdot l^{-1})$, $g\in GL_{n},$ which is of course
isomorphic to $\rho $. Now, the fact that $L_{\left\lfloor \frac{n}{2}%
\right\rfloor }$ normalizes $U=U_{\left\lfloor \frac{n}{2}\right\rfloor }$,
and acts transitively on each collection $\mathcal{O}_{r}$ of all matrices
in $U$ of a fixed rank, completes the verification of the proposition.
\end{proof}

\subsection{\textbf{Proofs for Section \protect\ref{S-BtoGP}}}

\subsubsection{\textbf{Proof of Claim \protect\ref{C-Slgk}\label{P-C-Slgk}}}

\begin{proof}
Parts (2), and (3), of the claim follow using a direct substitution of the
numerical data given by Formulas (\ref{CRs-SLn}), (\ref{Dim-SLn}), and (\ref%
{Card-SLn}).

Concerning Part (1). First we know that for each member $\pi \in (\widehat{SL%
}_{n})_{\otimes ,1}$, we have,\vspace{0.01in}

\begin{itemize}
\item $\dim (\pi )=q^{n-1}+o(...);\vspace{0.01in}$

and,

\item $\frac{\chi _{\rho }(T)}{\dim (\pi )}=\frac{1}{q}+o(...).$
\end{itemize}

So we can write%
\begin{eqnarray}
(S_{\ell ,g})_{1} &=&\dsum\limits_{\pi \in (\widehat{SL}_{n})_{\otimes
,1}}\dim (\pi )\left( \frac{\chi _{\pi }(T)}{\dim (\pi )}\right) ^{\ell
}\chi _{\pi }(g)  \label{Sl1} \\
&=&\left( \frac{1}{q}(\frac{1}{q})^{\ell -n}+o(...)\right) \dsum\limits_{\pi
\in (\widehat{SL}_{n})_{\otimes ,1}}\chi _{\pi }(g).  \notag
\end{eqnarray}%
Now, to understand the right-hand side factor in (\ref{Sl1}), recall that,
the irreps of $SL_{n}$ of tensor rank $k=1$, are given in Example \ref%
{Ex-U-rank-k=1}. In particular, they are members of the space $L^{2}(\mathbb{%
F}_{q}^{n})$ of complex valued functions on $\mathbb{F}_{q}^{n}$. In more
detail, each of them appears there with multiplicity one, in addition to the
trivial rep $\mathbf{1}$ which shows up twice, i.e., we have, 
\begin{equation*}
L^{2}(\mathbb{F}_{q}^{n})=2\cdot \mathbf{1}+\dsum\limits_{\rho \in (\widehat{%
SL}_{n})_{\otimes ,1}}\rho .
\end{equation*}%
Recall that we denoted the permutation representation of $GL_{n}$ on that
space by $\omega _{n,1}$. In particular, we conclude that%
\begin{equation}
\dsum\limits_{\rho \in (\widehat{SL}_{n})_{\otimes ,1}}\chi _{\rho
}(g)=trace(\omega _{n,1}(g)\curvearrowright L^{2}(\mathbb{F}_{q}^{n}))-2=-2,
\label{-2}
\end{equation}%
taking into account that the element $g$ has no eigenvalue equal to $1$.

Combining (\ref{Sl1}) and (\ref{-2}), we see that Part (1) holds true. This
completes the verification of the claim.
\end{proof}

\subsection{\textbf{Proofs for Appendix \protect\ref{A-FT}}}

\subsubsection{\textbf{Proof of Claim \protect\ref{C-F-FT-Mmnk}\label%
{P-C-F-FT-Mmnk}}}

\begin{proof}
Each element of the set $(M_{m,n})_{k}^{o}(B)$ contributes $1$ to the FT (%
\ref{FT-Mmnk}) of $1_{(M_{m,n})_{k}}$ evaluated at $B$. On the other hand
the multiplicative group $\mathbb{F}_{q}^{\ast }$ acts naturally by scaling
on the complement $(M_{m,n})_{k}\smallsetminus (M_{m,n})_{k}^{o}(B)$, and
(by one-dimensional harmonic analysis, i.e., orthogonality of characters of $%
\mathbb{F}_{q}$) each orbit contributes $-1$ to the FT. Overall we have,%
\begin{eqnarray}
\widehat{1}_{(M_{m,n})_{k}}(B) &=&-(\frac{1}{q-1})(\#((M_{m,n})_{k})-%
\#((M_{m,n})_{k}^{o}(B)))+\#((M_{m,n})_{k}^{o}(B))  \label{F-FT-Mmnk-B} \\
&=&-(\frac{1}{q-1})\#((M_{m,n})_{k})+(\frac{q}{q-1})\#((M_{m,n})_{k}^{o}(B)),
\notag
\end{eqnarray}%
as claimed.
\end{proof}

\subsubsection{\textbf{Proof of Theorem \protect\ref{T-FT-Mmnk-T}\label%
{P-P-FT-Mmnk-T}}}

We want to use Formula (\ref{F-FT-Mmnk-B}) above. For the rank one matrix $T$
we denote $(M_{m,n})_{k}^{o}=(M_{m,n})_{k}^{o}(T)$, and we note that,

\begin{lemma}
\label{L-M,Mo}We have,

\begin{enumerate}
\item \label{L-M,Mo-P1}$\#((M_{m,n})_{k})=\#(\Gamma _{n,k})\#(\Gamma
_{m,k})\#(GL_{k});\smallskip $

and,

\item \label{L-M,Mo-P2}$\#((M_{m,n})_{k}^{o})=\left(
(q^{m-1}-1)(q^{n}-1)+(q^{n-k}-1)(q-1)q^{m-1}\right) (\frac{\#(\Gamma
_{n,k})\#(\Gamma _{m,k})\#(GL_{k})}{(q^{n}-1)(q^{m}-1)}).$
\end{enumerate}
\end{lemma}

We will verify Lemma \ref{L-M,Mo} below.

\begin{proof}
(of Theorem \ref{T-FT-Mmnk-T}) \ We compute%
\begin{eqnarray*}
\widehat{1}_{(M_{m,n})_{k}}(T) &=&-(\frac{1}{q-1})\#((M_{m,n})_{k})+(\frac{q%
}{q-1})\#((M_{m,n})_{k}^{o}) \\
&=&\left( q^{n+m-k}-q^{n}-q^{m}+1\right) (\frac{\#(\Gamma _{n,k})\#(\Gamma
_{m,k})\#(GL_{k})}{(q^{n}-1)(q^{m}-1)}),
\end{eqnarray*}%
where in the second equality we used Lemma \ref{L-M,Mo}.
\end{proof}

\begin{proof}
(of Lemma \ref{L-M,Mo}) \smallskip

\textit{Part \ref{L-M,Mo-P1}}. This is a consequence of the rank-nullity
theorem. Indeed, for $A$ in $(M_{m,n})_{k}$, the kernel is an $n-k$
dimensional subspace of $\mathbb{F}_{q}^{n}$, while the image is a $k$%
-dimensional subspace of $\mathbb{F}_{q}^{m}$, and $A$ defines an invertible
transformation from $\mathbb{F}_{q}^{n}/\ker (A)$ to $im(A)$.\smallskip\ The
assertion follows.

\textit{Part \ref{L-M,Mo-P2}. }First lets give some standard formula for a
general matrix of rank one in $M_{n,m}$. Fix a vector $0\neq v\in \mathbb{F}%
_{q}^{n}$, and a linear functional $0\neq \lambda $ in the dual space $(%
\mathbb{F}_{q}^{m})^{\ast }$. Then 
\begin{equation*}
T_{v,\lambda }:\mathbb{F}_{q}^{m}\rightarrow \mathbb{F}_{q}^{n},\text{ \ }%
T_{v,\lambda }(w)=\lambda (w)v,
\end{equation*}%
is a rank one operator, and any rank one matrix in $M_{n,m}$ is of this form.

We are trying to compute the cardinality of $(M_{m,n})_{k}^{o}=\{A\in
(M_{m,n})_{k};$ $trace(T_{v,\lambda }\circ A)=\lambda (A(v))=0\}.$

We will count and find (justification below) that the number of the matrices 
$A\in (M_{m,n})_{k}$ that,\smallskip

(a) satisfy $A(v)=0$, is $\#(\Gamma _{n-1,k})\#(\Gamma
_{m,k})\#(GL_{k});\smallskip $

(b) satisfy $im(A)\subset \ker (\lambda ),$ is $\#(\Gamma _{n,k})\#(\Gamma
_{m-1,k})\#(GL_{k});\smallskip $

(c) satisfy both (a) and (b) above, is $\#(\Gamma _{n-1,k})\#(\Gamma
_{m-1,k})\#(GL_{k});\smallskip $

(d) satisfy $\lambda (A(v))=0,$ but neither (a) nor (b), is $(\#(\Gamma
_{n,k})-\#(\Gamma _{n-1,k}))((\frac{q^{m-1}-1}{q^{m}-1})\#(\Gamma
_{m,k})-\#(\Gamma _{m-1,k}))\#(GL_{k}).$

In addition, Formula (\ref{Card-Gama-nk}) implies that 
\begin{equation}
\#(\Gamma _{n-1,k})=(\frac{q^{n-k}-1}{q^{n}-1})\#(\Gamma _{n,k}).
\label{Gama_nk-Gama_n-1k}
\end{equation}%
Now, note that 
\begin{equation*}
\#((M_{m,n})_{k}^{o})=\text{(a)+(b)-(c)+(d),}
\end{equation*}%
and so, a direct calculation using the explicit cardinalities presented in
(a)-to-(d), including Identity (\ref{Gama_nk-Gama_n-1k}), produce the
assertion made in Part \ref{L-M,Mo-P2}.\smallskip

Let us now finish the proof, by justifying (a)-to-(d).

We will frequently use the fact that $\#(\Gamma _{n,n-k})=\#(\Gamma _{n,k})$%
.\smallskip

To justify (a), note that if $A(v)=0$, then $\ker (A)$ defines $(n-k-1)$%
-dimensional subspace of $\mathbb{F}_{q}^{n}/span(v)\simeq \mathbb{F}%
_{q}^{n-1}$, and there are $\#(\Gamma _{n-1,n-k-1})=\#(\Gamma _{n-1,k})$
such subspaces. On the other hand, $im(A)$ is a $k$-dimensional subspace of $%
\mathbb{F}_{q}^{n}$, and each of these $\#(\Gamma _{m,k})$ is a legitimate
one. After these choices were made, $A$ can be viewed as an element of $%
GL_{k}$, and again each such is possible. Overall, Part (a)
follows.\smallskip

The derivation of Parts (b) and (c) is very similar to that of Part (a), so
let us omit it.\smallskip

Finally, let us justify Part (d). For the $A$'s there the kernel is $n-k$
dimensional subspace, so a member of $\Gamma _{n,n-k}$, and there are $%
\#(\Gamma _{n,n-k})=\#(\Gamma _{n,k})$ such. But $v$ is not allowed to be in 
$\ker (A)$, for such $A$, so is not one of the subspaces that are sent to an 
$n-k-1$ dimensional subspace of $\mathbb{F}_{q}^{n}/span(v)\simeq \mathbb{F}%
_{q}^{n-1}$, and there are $\#(\Gamma _{n-1,n-k-1})=\#(\Gamma _{n-1,k})$ of
such. After one is making one of these $\Gamma _{n,k}\smallsetminus \Gamma
_{n-1,k}$ choices, we should look on the constraints on the image of $A$.
The first one is that $A(v)\in \ker (\lambda )\smallsetminus \{0\}$, and we
have $q^{m-1}-1$ such options. After this choice, we need to make sure $A$
has rank $k$. For this the number of options is%
\begin{eqnarray*}
(q^{m}-q)\cdot ...\cdot (q^{m}-q^{k-1}) &=&(\frac{1}{q^{m}-1}%
)\#((M_{m,n})_{k}) \\
&=&(\frac{1}{q^{m}-1})\#(\Gamma _{m,k})\#(GL_{k}).
\end{eqnarray*}%
And, lastly, we need to make sure that $im(A)$ is not contained in $\ker
(\lambda )$, so we need to subtract $\#(\Gamma _{m-1,k})$ options.

Overall, when we multiply the number of options for the domain of $A$ with
that for its range, we get what is claimed in Part (d). This completes the
proof of Lemma \ref{L-M,Mo}.
\end{proof}


\begin{thebibliography}{Bezrukavnikov-Liebec-Shalev-Tiep18}
\bibitem[Arad-Herzog-Stavi85]{Arad-Herzog-Stavi85} Arad Z., Herzog M., and
Stavi J., Powers and products of conjugacy classes in groups. \textit{\href{https://link.springer.com/chapter/10.1007/BFb0072286}%
{ LNM 1112, Springer-Verlag, Berlin, (1985) 6-51}}.

\bibitem[Artin57]{Artin57} Artin E., Geometric Algebra. \href{https://archive.org/details/geometricalgebra033556mbp}%
{\textit{Interscience, New York (1957)}}.

\bibitem[Auslander-Tolimieri79]{Auslander-Tolimieri79} Auslander L. and
Tolimieri R., Is computing with the finite Fourier transform pure or applied
mathematics? \textit{\href{https://projecteuclid.org/journals/bulletin-of-the-american-mathematical-society-new-series/volume-1/issue-6/Is-computing-with-the-finite-Fourier-transform-pure-or-applied/bams/1183544898.full}%
{Bulletin of the AMS Vol. 5 (1981), 263-312.}}

\bibitem[Bezrukavnikov-Liebec-Shalev-Tiep18]%
{Bezrukavnikov-Liebec-Shalev-Tiep18} Bezrukavnikov R., Liebeck M., Shalev
A., and P. H. Tiep. Character bounds for finite groups of Lie type. \href{https://www.intlpress.com/site/pub/pages/journals/items/acta/content/vols/0221/0001/a001/index.php}%
{\textit{Acta Math. 221 (2018) 1-57}}.

\bibitem[Borel69]{Borel69} Borel A., Linear algebraic groups. \href{https://link.springer.com/book/10.1007/978-1-4612-0941-6}%
{\textit{GTM 126, Springer-Verlag (1969)}}\textit{.}

\bibitem[Bump04]{Bump04} Bump D., Lie Groups. \href{https://www.springer.com/gp/book/9781475740943}%
{\textit{Springer, New York (2004).}}

\bibitem[Clifford37]{Clifford37} Clifford A. H., Representations induced in
an invariant subgroup. \href{https://www.jstor.org/stable/1968599?seq=1}{%
\textit{Annals of Math 38 (1937) 533--550.}}

\bibitem[Deligne-Lusztig76]{Deligne-Lusztig76} Deligne P. and Lusztig G., 
\textit{\href{https://www.jstor.org/stable/1971021?seq=1}{Representations of
reductive groups over finite fields.\textit{\ Annals of Math. 103 (1976),
103-161}}.}

\bibitem[Frobenius1896]{Frobenius1896} Frobenius F.G., \"{U}ber
Gruppencharaktere. \textit{Sitzber. Preuss. Akad. Wiss. (1896) 985--1021.}

\bibitem[Fulton-Harris91]{Fulton-Harris91} Fulton W. and Harris J.,
Representation theory: A first course. \href{https://www.springer.com/gp/book/9780387975276}%
{\textit{GTM 129, Springer (1991).}}

\bibitem[Gel'fand70]{Gel'fand70} Gel'fand S.I., Representations of the full
linear group over a finite field. \href{https://iopscience.iop.org/article/10.1070/SM1970v012n01ABEH000907}%
{\textit{Math. USSR-Sb., 12 (1970) 13-39}}\textit{.}

\bibitem[Gerardin77]{Gerardin77} G\'{e}rardin P., Weil representations
associated to finite fields. \textit{\ \href{https://www.sciencedirect.com/science/article/pii/0021869377903945}%
{\textit{J\textit{. Alg. 46 (1977), 54-101}}}.}

\bibitem[Green55]{Green55} Green J.A., The characters of the finite general
linear groups. \href{https://www.ams.org/journals/tran/1955-080-02/S0002-9947-1955-0072878-2/home.html}%
{\textit{TAMS 80 (1955) 402--447}}\textit{.}

\bibitem[Gurevich-Howe15]{Gurevich-Howe15} Gurevich S. and Howe R., Small
Representations of finite classical groups. \textit{\href{https://link.springer.com/chapter/10.1007/978-3-319-59728-7_8}%
{\textit{Proceedings of Howe's 70th birthday conference, New Haven (2015)}}.}

\bibitem[Gurevich-Howe17]{Gurevich-Howe17} Gurevich S. and Howe R., Rank and
Duality in Representation Theory.\textit{\ \href{https://link.springer.com/article/10.1007/s11537-020-1728-3}%
{The 19th Takagi Lectures - July 2017, Jpn J. Math. 15 (2020) 223-309}}%
\textbf{.}

\bibitem[Gurevich-Howe18]{Gurevich-Howe18} Gurevich S. and Howe R., A look
on representations of $SL_{2}(\mathbb{F}_{q})$ through the lens of size. 
\href{https://link.springer.com/article/10.1007/s40863-018-0098-8}{\textit{%
Joe Wolf's 80th Birthday Volume, S\~{a}o Paulo J. of Math. Sci. 12 (2018)
252--277.}}

\bibitem[Gurevich-Howe19]{Gurevich-Howe19} Gurevich S. and Howe R., Harmonic
Analysis on $GL_{n}$\ over Finite Fields.\textit{\ \href{https://arxiv.org/abs/2105.12369}%
{\textit{The Bertram Kostant Memorial Volume. \textit{Pure and Applied
Mathematics Quarterly (2019) 62p. }Accepted.}}}

\bibitem[Harish-Chandra70]{Harish-Chandra70} Harish-Chandra., Eisenstein
series over finite fields. \textit{\href{https://link.springer.com/chapter/10.1007/978-3-642-48272-4_3}%
{\textit{Functional analysis and related fields, Springer (1970) 76--88}}. }

\bibitem[Harish-Chandra84]{Harish-Chandra84} Harish-Chandra., Collected
papers IV - 1970-83. \textit{\href{https://www.springer.com/gp/book/9783662454350}%
{Springer-Verlag (1984)}}.

\bibitem[Hartshorne77]{Hartshorne77} Hartshorne R., Algebraic Geometry. 
\textit{\href{https://www.springer.com/us/book/9780387902449}{%
Springer-Verlag, GTM 52 (1977).}}

\bibitem[Howe73]{Howe73} Howe R., Invariant theory and duality for classical
groups over finite fields with applications to their singular representation
theory. \textit{Preprint, Yale University (1973). }

\bibitem[Howe-Moy86]{Howe-Moy86} Howe R. and Moy A., Harish-Chandra
homomorphisms for p-adic groups. \href{https://bookstore.ams.org/cbms-59}{%
\textit{CBMS Regional Conference Series in Mathematics 59 (1986).}}

\bibitem[Humphries80]{Humphries80} Humphries S.P., Generation of special
linear groups by transvections. \href{https://www.sciencedirect.com/science/article/pii/0021869386900414}%
{\textit{J. of Alg. 99 (1986) 480-495}}\textit{.}

\bibitem[Lang-Weil54]{Lang-Weil54} Lang S. and Weil A. Number of points of
varieties in finite fields.\textit{\ \href{https://mathscinet.ams.org/mathscinet-getitem?mr=65218}%
{Amer. J. Math. 76 (1954), 819--827}.}

\bibitem[Lusztig84]{Lusztig84} Lusztig, G., Characters of Reductive Groups
over a Finite Field. \href{https://www.jstor.org/stable/j.ctt1b9x10c}{%
\textit{Annals of Math. Studies, Princeton University Press (1984)}}\textit{.%
}

\bibitem[Mackey49]{Mackey49} Mackey G.W., Imprimitivity for representations
of locally compact groups I. \href{https://www.jstor.org/stable/1969423?seq=1}%
{\textit{PNAS 35 (1949) 537--545.}}

\bibitem[Schur1905]{Schur1905} Schur I., Neue Begr\"{u}ndung der Theorie der
Gruppencharaktere. \textit{Sitzungsberichte der K\"{o}niglich Preu\ss ischen
Akademie der Wissenschaften zu Berlin (1905) 406-432. }

\bibitem[Serre77]{Serre77} Serre J.P., Linear Representations of Finite
Groups. \href{https://link.springer.com/book/10.1007%2F978-1-4684-9458-7}{%
\textit{Springer (1977)}}\textit{.}

\bibitem[Weil64]{Weil64} Weil A., Sur certains groups d'operateurs
unitaires. \href{https://link.springer.com/article/10.1007/BF02391012}{%
\textit{Acta Math. 111 (1964), 143-211.}}

\bibitem[Zelevinsky81]{Zelevinsky81} Zelevinsky A., Representations of
finite classical groups: A Hopf algebra approach. \href{https://www.springer.com/gp/book/9783540108245}%
{\textit{LNM 869. Springer (1981).}}
\end{thebibliography}
\end{document}